\documentclass{amsart}

\usepackage{amssymb}
\usepackage[dvips]{epsfig}
\usepackage{graphicx}
\usepackage{color}

\newcommand{\R}{{\mathbb R}}

\newtheorem{theorem}{Theorem}[section]

\theoremstyle{definition}

\newtheorem*{merci}{Acknowledgements}

\theoremstyle{remark}
\newtheorem{remark}{Remark}[section]
\def\R{{\mathbb R}}

\numberwithin{equation}{section}

\begin{document}
\title[On the Kaup-Broer-Kupershmidt  systems]{On the Kaup-Broer-Kupershmidt  systems }

\author[C. Klein]{Christian Klein}
\address{Institut de Math\'ematiques de Bourgogne,  UMR 5584;\\
Institut Universitaire de France \\
Universit\'e de Bourgogne, 9 avenue Alain Savary, 21078 Dijon
                Cedex, France} %\\
\email{Christian.Klein@u-bourgogne.fr}

\author[J.-C. Saut]{Jean-Claude Saut}
\address{Laboratoire de Math\' ematiques, UMR 8628,\\
Universit\' e Paris-Saclay et CNRS\\ 91405 Orsay, France}
\email{jean-claude.saut@universite-paris-saclay.fr}

      \date{February 23rd, 2024}                                                                                                                                                                                                      

\maketitle
\begin{center}
In memoriam Thomas Kappeler (1953-2022)

\end{center}

\begin{abstract}
%\textit{Abstract}
The aim of this paper is to survey and complete, mostly by numerical simulations, results on a  remarkable Boussinesq system  describing weakly nonlinear, long surface water waves. It is the only member of  the so-called (abcd) family of Boussinesq systems  known to be completely integrable.

\end{abstract}

\section{introduction}

We are interested in this paper in a particular case of the so-called 
abcd Boussinesq systems for surface water waves, see \cite{BCL, 
BCS1,BCS2} \footnote{Boussinesq \cite{Bou2} was the first to derive a particular Boussinesq system, not in the class of those studied here though. We refer to \cite{Da,Da2,Da3} for details and for an excellent  history of hydrodynamics in the nineteenth century.}
\begin{equation}
    \label{abcd2}
    \left\lbrace
    \begin{array}{l}
    \eta_t+\nabla \dot {\bf v}+\epsilon \nabla\cdot(\eta {\bf 
	v})+\mu\lbrack a \nabla\cdot \Delta{\bf v}-b\Delta \eta_t\rbrack=0, \\
    {\bf v}_t+\nabla \eta+\epsilon \frac{1}{2}\nabla |{\bf v}|^2+\mu\lbrack c\nabla \Delta \eta-d\Delta {\bf v}_t\rbrack=0.
\end{array}\right.
    \end{equation}
Here $\eta= \eta(x,t), x\in \R^d, d=1,2, t\in \R$ is the elevation of the wave, ${\bf v}={\bf v}(x,t)$ is a measure of the horizontal velocity, $\mu$ and $\epsilon$ are the small parameters (shallowness and nonlinearity parameters respectively) defined as
$$\mu=\frac{h^2}{\lambda^2}, \quad \epsilon= \frac{\alpha}{h}$$
where  $\alpha$ is a typical amplitude of the wave, $h$ a typical depth and $\lambda$ a typical horizontal wavelength.

In the Boussinesq regime, $\epsilon$ and $\mu$ are supposed to be of 
same order, $\epsilon\sim\mu\ll1,$ and we will take for simplicity $\epsilon=\mu,$  writing \eqref{abcd2} as \footnote{Particular cases are formally derived in \cite{Broer, Ding, Per}.}
\begin{equation}
    \label{abcd}
    \left\lbrace
    \begin{array}{l}
    \eta_t+\nabla \cdot {\bf v}+\epsilon \lbrack\nabla\cdot(\eta {\bf 
	v})+a \nabla\cdot \Delta{\bf v}-b\Delta \eta_t\rbrack=0, \\
    {\bf v}_t+\nabla \eta+\epsilon\lbrack \frac{1}{2}\nabla |{\bf v}|^2+c\nabla \Delta \eta-d\Delta {\bf v }_t\rbrack=0,
\end{array}\right.
    \end{equation}

The coefficients (a, b, c, d) are restricted by the condition 
$$a+b+c+d=\frac{1}{3}-\tau,$$
where $\tau\geq 0$ is the surface tension coefficient.

When restricted to one-dimensional, unidirectional motions, 
\eqref{abcd} leads to the Korteweg-de Vries (KdV) equation, see \cite{La1}

$$u_t+u_x+\epsilon\left(\frac{1}{3}-\tau\right)u_{xxx}+\epsilon uu_x=0.$$

The class of systems \eqref{abcd2}, \eqref{abcd} models  water waves on a flat bottom propagating in both directions in the aforementioned  regime (see \cite{BCL, BCS1, BCS2}).

It turns out that two particular one-dimensional cases of the abcd 
systems have remarkable properties. The first one, we will refer to 
as the {\it Amick-Schonbek system}, can be viewed as a dispersive perturbation of the Saint-Venant (shallow water) system\footnote{ Actually this system is a particular case of a system derived by Peregrine in \cite{Per} but Schonbek and Amick were the first to point out its remarkable mathematical properties.}. We refer to \cite{KS-AS} for a review of known results together with new numerical simulations.

 The second one, referred to as the {\it Kaup-Broer-Kupershmidt system} (KBK)  was introduced in \cite {Broer, BGT}, as    a   surface long wave model (see   equations (3.2), (3.3) in \cite{Broer}) and recently in \cite{CFOPV} as an internal wave model. On the other hand Kaup and Kupershmidt introduced it as an integrable system \cite{Ka, Kup}.  It corresponds in the abcd family to $a=\pm 1, b=c=d=0$ and writes

\begin{equation}
    \label{KBK}
    \left\lbrace
    \begin{array}{l}
    \eta_t+v_x+(\eta v)_x+\alpha v_{xxx}=0, \\
    v_t+\eta_x+vv_x=0.
\end{array}\right.
    \end{equation}

    $\alpha=1$ corresponds to the {\it bad} KBK system and $\alpha 
	=-1$ to the {\it good} KBK system. The two systems turn out to be 
	integrable in a sense that will be detailed below.

    \begin{remark}
    When viewed as water wave model, the variable $\eta$ in \eqref{KBK} represents the elevation of the wave, so that physically the total depth $\zeta=1+\eta$ should be positive. The non-cavitation  condition $\zeta>0$  will be mostly  ignored in what follows. Note that in terms of $(v,\zeta)$ \eqref{KBK} writes
    \begin{equation}
    \label{KBKbis}
    \left\lbrace
    \begin{array}{l}
    \zeta_t+(\zeta v)_x+\alpha v_{xxx}=0, \\
    v_t+\zeta_x+vv_x=0.
\end{array}\right.
    \end{equation}
    
    Setting $U=v_x$ in $\eqref{KBKbis}_1$ and solving for $\zeta$ leads to the {\it Boussinesq like} equation
    
   \begin{equation}\label{Bousslike}
   U_{tt}+U_tU_{xx}+2U_xU_{xt}+\frac{3}{2}U_x^2U_{xx}-\alpha  U_{xxxx}=0.
   \end{equation} 
    
    \end{remark}

    \begin{remark}
    The KBK system is sometimes written in terms of the velocity potential $\phi$ such that $\phi_x(x,t)=v(x,t)$:
    
    \begin{equation}
    \label{KBK2}
    \left\lbrace
    \begin{array}{l}
    \eta_t+\phi_{xx}+(\eta \phi_x)_x+\alpha \phi_{xxxx}=0, \\
    \phi_t+\eta+\frac{1}{2}(\phi_x)^2=0.
\end{array}\right.
    \end{equation}

    \end{remark}
    \begin{remark}
    The KBK system should not be confused with the so-called 
     Kaup-Kupershmidt equation
    \begin{equation}\label{KK-eq}
     v_t=v_{xxxxx}+10vv_{xxx}+25v_xv_{xx}+20v^2v_x,
    \end{equation}
    that was introduced in \cite{Kaup2, Kup2}, and which is the first 
	equation in a hierarchy (different from the KdV hierarchy) of 
	integrable equations with Lax operator 
	$\partial_x^3+2u\partial_x+u_x$, see for instance \cite{GWH}. We 
	are not aware of a physical application of the Kaup-Kupershmidt equation.
    \end{remark}
\vspace{0.3cm}

The paper is organized as follows. In   section 2 we will describe 
general facts on the   KBK system.  The next section reviews the 
results obtained by partial differential equations (PDE) techniques 
while section 4 focusses on the integrability side. In section 5 we 
introduce a numerical approach to the KBK system and test it for the 
known soliton. In section 6 the known stability of the solitons is 
illustrated. In section 7 the long time behavior of  KBK solutions 
for localised initial data is studied showing that the soliton 
resolution conjecture can be applied. In section 8 we explore the 
formation of dispersive shock waves in the vicinity of shocks of the 
corresponding dispersionless Saint-Venant system. We add some 
concluding remarks in section 9.

\section{Generalities on the Kaup-Broer-Kupershmidt system}

We recall here that the Kaup-Broer-Kupershmidt system is the one-dimensional version of the 2-D (abcd Boussinesq) system
\begin{equation}
    \label{Kaup-true-2D}
    \left\lbrace
    \begin{array}{l}
    \eta_t+\nabla \cdot {\bf v}+\epsilon \nabla\cdot(\eta {\bf v})+\alpha\epsilon \Delta \nabla\cdot {\bf v}=0,\; \alpha=\pm 1 \\
    {\bf v}_t+\nabla \eta+ \frac{\epsilon}{2}\nabla |{\bf v}|^2=0,
\end{array}\right.
    \end{equation}

which is linearly well-posed when $\alpha=-1.$ The local well-posedness of the Cauchy problem is established in \cite{KS} (see also \cite {SWX}). The more difficult question of  {\it long time existence} (that is on time scales of order $ O(1/\epsilon))$ is established in \cite{SWX} under the non-cavitation condition $1+\epsilon \eta>0.$ 

On the other hand \eqref{Kaup-true-2D} is ill-posed when $\alpha=+1.$ Note that the ill-posed version can be turned into  a well-posed system by using the BBM trick, leading to

\begin{equation}
    \label{Kaup-true-mod}
    \left\lbrace
    \begin{array}{l}
    \eta_t+\nabla \cdot {\bf v}+\epsilon \nabla\cdot(\eta {\bf v})-\epsilon\Delta \eta_t=0 \\
    {\bf v}_t+\nabla \eta+\epsilon \frac{1}{2}\nabla |{\bf v}|^2=0,
\end{array}\right.
    \end{equation}
    for which the long time existence  is established in \cite{SWX}.
    
    \vspace{0.5cm}
    
    We will from now on restrict to the one-dimensional case with most of the time  $\epsilon =1$. The linearized  "bad" KBK system ($\alpha=1$ in 
    \eqref{KBK})  is Hadamard ill-posed since  the dispersion relation is $\omega^2=k^2- k^4$, showing the short wave ill-posedness. It was established in \cite{ABM} that nonlinearity does not erase the problem, making the Cauchy problem ill-posed in all Sobolev spaces in the sense that  arbitrary small smooth solutions can blow-up in arbitrary short time in Sobolev spaces norms.
    
     \begin{remark}
 Following \cite{Kup}, the change of variable $v=v, \eta=h-v_x$ transforms the Kaup-Broer-Kupershmidt system with $\alpha=1$  into the (ill-posed) system
  
 \begin{equation}
    \label{new-BKK}
    \left\lbrace
    \begin{array}{l}
    h_t+v_x+ h_{xx}+(vh)_x=0, \\
    v_t+h_x- v_{xx}+ vv_x=0
    \end{array}\right.
    \end{equation}
with dispersion $\omega(k)=\pm i k\sqrt{1-k^2}.$
   
On the other hand, in the well-posed case $\alpha=-1,$ the change of variable $v=v, \eta=h+iv_x$ leads to the (linearly well-posed)  "Schr\"{o}dinger like" system

 \begin{equation}
    \label{new-BKK2}
    \left\lbrace
    \begin{array}{l}
    h_t+v_x-i h_{xx}+(vh)_x=0, \\
    v_t+h_x+i v_{xx}+ vv_x=0,
    \end{array}\right.
    \end{equation}
    with dispersion $\omega(k)=\pm ik(1+k^2)^{1/2}$ 

 \end{remark}

As recalled in \cite{NZ},   Broer \cite{Broer} and Kaup \cite{Ka} derived the system in dimensional form
\begin{equation}
    \label{BKK}
    \left\lbrace
    \begin{array}{l}
    \eta_t+h_0\phi_{xx}+(\eta\phi_x)_x+\left(\frac{h_0^3}{3}-\frac{h_0\tau}{\rho g}\right)\phi_{xxx}=0 \\
    \phi_t+\frac{1}{2}(\phi_x)^2+g\eta=0,
\end{array}\right.
    \end{equation}
where $\eta$ is the elevation of the wave, $\phi$ the velocity potential evaluated at the free surface, the constant $h_0$ is the quiescent water depth, $g>0$  is the gravitational acceleration, $\tau\geq 0$ is the surface tension coefficient and $\rho$ the density of the fluid.

The "bad" version corresponds to pure gravity waves ($\tau=0$) or to 
gravity-capillarity waves with small surface tension 
($\tau/(g\rho)<h_0^2/3$, corresponding to $b=d=c=0, a>0$ in \eqref{abcd}). The "good" version occurs with strong surface tension that is when $\tau/(g\rho)>h_0^2/3,$ (a not too physical case...), corresponding to $b=d=c=0, a<0$ in \eqref{abcd} . The difference of nature of the equation depending on the capillarity  is reminiscent of that of the Kadomtsev-Petviashvili equation (KP-II for gravity waves and KP-I for gravity-capillarity waves with strong surface tension) although both KP-II and KP-I are well-posed!

 \begin{remark}
 Despite the fact that the Kaup-Broer-Kupershmidt system with the 
 plus sign  is linearly ill-posed, it possesses solitary wave solutions $(\eta(x-kt), v(x-kt))$ when $k>1$ \footnote{In the classical version $a=1/3$.}.  This has been found by Kaup, \cite{Ka} using the Inverse Scattering machinery who also proved the existence of N-solitons. Matveev and Yavor \cite{MY}, exhibited a vast class of almost periodic solutions containing N-soliton solutions  as a degenerate case, see also \cite{NZ}. The existence of undular bores is proven in \cite{EGP}.
 
 A direct approach to the existence of solitary waves  is given in \cite{chen} with $a=\frac{1}{3}$. Actually a solitary wave satisfies the equations
 
 $$(v')^2=R_1(u)\equiv \frac{3}{4}v^2(v-(2k-2))(v-(2k+2)) ,\quad \eta=v\left(k-\frac{v}{2}\right),$$
where the derivative is taken with respect to $\xi=x-kt.$

By studying the function $R_1,$ one shows that there exists a unique solution $v(\xi)$ which is even and monotone for $\xi>0$ for any $k>1$. More precisely, the solution $v(\xi)$ 
reads
$$
v(\xi)=\frac{2(k^2-1)}{\cosh\bigl(\sqrt{3(k^2-1)}\xi\bigr)+k}.
$$
We also have
$$
\eta(\xi)=\frac{2(k^2-1)\Bigl(k\cosh\bigl(\sqrt{3(k^2-1)}\xi\bigr)+1\Bigr)}{\Bigl(\cosh\bigl(\sqrt{3(k^2-1)}\xi\bigr)+k\Bigr)^2}.
$$

When $1<k\leq 2,$ the corresponding $\eta(\xi)$ is also monotonically decreasing for $\xi>0$ while this is no longer true when $k>2.$ 

One notices that for $1<k\leq 2,$ $\|v\|_{L^\infty}=\|\eta\|_{L^\infty}=2(k-1).$
 
 \vspace{0.3cm}
 \end{remark}

\section{The Kaup-Broer-Kupershmidt system by PDE techniques}

We recall that the Cauchy problem  for \eqref{Kaup-true-2D}  with $\alpha=-1$ is well-posed  in both spatial dimensions one and two under the non-cavitation condition $1+\epsilon \eta>0$ on time scales of order $O(1/\epsilon),$ see \cite{SWX}.

We now focus on the one-dimensional "good" KBK system with $\epsilon =1$ which has the Hamiltonian structure

$$\eta_t+\partial_x\frac{\delta \mathcal H}{\delta v}=0,\quad v_t+\partial_x\frac{\delta \mathcal H}{\delta \eta}=0,$$

where the Hamiltonian $\mathcal H$ is given by 
$$
\mathcal H=\frac{1}{2}\int [\eta^2+(1+\eta)v^2+(v_x^2)]dx.$$

\begin{remark}
The Hamiltonian of the "bad" KBK system is
$$
\mathcal H=\frac{1}{2}\int [\eta^2+(1+\eta)v^2-(v_x^2)]dx.$$

\end{remark}

\begin{remark}
Both the "bad" and "good" KBK system have the conserved quantity:
$$H_0=\int_\R\eta v dx.$$
\end{remark}

Angulo, \cite{Ang}, proved that the Cauchy problem for the "good" KBK system  is locally well-posed in $H^{s-1}(\R)\times H^s(\R), s>3/2.$\footnote{A local well-posedness result in a weighted space with less regularity was obtained in \cite{KS}.} 

Surprisingly he also found  an additional conservation law namely

\begin{equation}\label{Ang}
\begin{split}
I_3(\eta,v)=&\frac{1}{8}\int [4(v_{xx}^2)+8(v_x^2)+4v^2+4(\eta_x^2)+4\eta^2+6v^2(v_x)^2 \\
&-16\eta vv_{xx}-4\eta(v_x)^2+10\eta v^2+2\eta^3+v^4+6\eta^2v^2+\eta v^4]dx.
\end{split}
\end{equation}

The existence of such a conservation law suggests that the "good" 
Kaup-Broer-Kupershmidt system  is integrable in the sense of Inverse 
Scattering. We refer to the next section for further information on this issue.

Actually $I_3$ can be viewed as the Hamiltonian of the (linearly well-posed)  higher order  system

\begin{equation}
    \label{KBKJaime}
    \left\lbrace
    \begin{array}{l}
    v_t-\eta_{xxx}+\eta_x+\lbrace\frac{5}{4}v^2+\frac{3}{4}\eta^2+\frac{3}{2}\eta v^2+\frac{1}{8}v^4-2vv_{xx}-\frac{1}{2}v_x^2\rbrace_x=0, \\
    \eta_t+v_{xxxxx}-2v_{xxx}+v_x+\lbrace\frac{3}{2}vv^2_x-\frac{3}{2}(v^2v_x)_x-2\eta v_{xx}-2(\eta v)_{xx}+\frac{1}{2}(\eta v_x)_x+\frac{5}{2}v\eta\\+\frac{1}{2}v^3+\frac{3}{2}(v\eta^2)+\frac{1}{2}v^3\eta\rbrace_x=0
\end{array}\right.
 \end{equation}
 
 that can be written
$$\eta_t+\partial_x\frac{\delta I_3}{\delta v}=0,\quad v_t+\partial_x\frac{\delta I_3}{\delta \eta}=0.$$

We are not aware of results on the Cauchy problem for \eqref{KBKJaime}.

  \vspace{0.3cm}
  In order to get a further insight into the "good" KBK system, it is useful to diagonalize the linear part to obtain the equivalent system, see \cite{SWX}:
  
\begin{equation}
    \label{bis}
    \left\lbrace
    \begin{array}{l}
    \eta_t+J_\epsilon \zeta_x+\frac{\epsilon}{2} N_1(\eta,w)=0, \\
    w_t-J_\epsilon w_x+\frac{\epsilon}{2} N_2(\eta,w)=0.
    \end{array}\right.
    \end{equation}
where $J_\epsilon =(I-\epsilon\partial_x^2)^{1/2}$ and

$$N_1(\eta,w)=\partial_x\lbrack 
(\zeta+w)J_\epsilon^{-1}(\eta-w)\rbrack+J_\epsilon\lbrack J_\epsilon^{-1}(\eta-w)J_\epsilon^{-1}(\eta_x-w_x)\rbrack ,$$
and
$$ N_2(\eta,w)= \partial_x\lbrack 
(\eta+w)J_\epsilon^{-1}(\zeta-w)\rbrack-J_\epsilon\lbrack J_\epsilon^{-1}(\eta-w)J_\epsilon^{-1}(\eta_x-w_x)\rbrack.  $$

Since
$$(1+\epsilon \xi^2)^{1/2}-\epsilon^{1/2}|\xi|=\frac{1}{(1+\epsilon\xi^2)^{1/2}+\epsilon^{1/2}|\xi|},$$
\eqref{bis} writes
\begin{equation}
    \label{ter}
    \left\lbrace
    \begin{array}{l}
    \eta_t+\epsilon^{1/2}\mathcal H\eta_{xx}+R_\epsilon \eta+\frac{\epsilon}{2}N_1(\eta,w)=0, \\
    w_t-\epsilon^{1/2}\mathcal H w_{xx}-R_\epsilon w+\frac{\epsilon}{2}N_2(\eta,w)=0.
    \end{array}\right.
    \end{equation}
where $R_\epsilon$ is the  (order zero) skew-adjoint  operator with symbol $\frac{i\xi}{(1+\epsilon\xi^2)^{1/2}+\epsilon^{1/2}|\xi|}.$

 Thus the "good" KBK system is equivalent to a system having a decoupled "Benjamin-Ono type" linear part. This property allows to apply the techniques used to obtain the local well-posedness of the Benjamin-Ono equation in low regularity spaces. Note that this approach would yield existence on the "short" time range $O(1/\sqrt \epsilon).$ See the discussion in \cite{SWX}.

  \begin{remark}
 In a recent paper, Melinand \cite{Mel} derived various dispersive 
 estimates  (Strichartz, local Kato smoothing, Morawetz) for a large class of  (abcd) systems including the Kaup-Broer-Kupershmidt system, both in one and two spatial dimensions. These estimates should play an important role for the local resolution of the Cauchy problem in "large" functional spaces and in the proof of the possible scattering of small solutions.
 \end{remark}
 
  We conclude this section by alluding to results in \cite{CIKP}, concerning the evolution of  discontinuous initial data under the flow of  the good KBK system.

   Actually a Riemann-type problem is considered with initial data 
$$\eta(x,t=0)=\eta_L,\quad\text{and}\quad v(x,t=0)=v_L\quad\text{for}\quad x<0,$$
$$\eta(x,t=0)=\eta_R,\quad\text{and}\quad v(x,t=0)=v_R\quad\text{for}\quad x>0,$$
and one finds in \cite{CIKP} a classification of the wave patterns 
evolving from these initial discontinuities.

\subsection{Traveling wave solutions}
   Following Angulo \cite{Ang} who considered the case $\epsilon =1,$ we now consider the traveling wave solutions of \eqref{KBK} that is solutions of the form ${\bf U}_c=(\eta,v)$ with
$$\eta(x,t)=n(x-ct),\; v(x,t)=\phi(x-ct)$$
with $c\in \mathbb{R}$ constant
so that $(n,\phi)$ satisfies the system (we have assumed that $(n,\phi)$ vanish at infinity which is natural in the water waves context):
\begin{equation}
   \left\lbrace
    \begin{array}{l}
    -cn+\phi+\epsilon n\phi-\epsilon\phi''=0, \\
    -c\phi+n+\epsilon \frac{\phi^2}{2}=0,
    \end{array}\right.
    \end{equation}
where $'$ denotes the derivative with respective to $\xi=x-ct$.

Eliminating $n$ from the second equation yields
\begin{equation}\label{SWSc1}
-\epsilon \phi''+(1-c^2)\phi+\frac{3}{2}C\epsilon \phi^2-\frac{\epsilon^2}{2}\phi^3=0.
\end{equation}

From Theorem 5 in \cite{BeLi} non trivial solitary waves exist if and 
only if $|c|<1.$ Actually \eqref{SWSc1}  can be integrated by 
standard methods, leading, for $|c|<1,$  a unique even solution (up to translations) which reads when $\epsilon =1$, see \cite{Ang}:
\begin{equation}\label{SWexpl}
v_{c,1}(\xi)=\frac{2(1-c^2)}{\cosh (\sqrt{1-c^2}\xi)-c},\quad \eta_{c,1}(\xi)=cv_{c,1}(\xi)-\frac{1}{2}v_{c,1}(\xi)
\end{equation}
from which we deduce the expression of $v_{c,\epsilon}:$
\begin{equation}\label{SWexp2}
v_{c,\epsilon}(\xi)=\frac{2(1-c^2)}{\epsilon[{\cosh 
(\sqrt{1-c^2}\epsilon^{-1/2}\xi)}-c]}
\end{equation}

   \vspace{0.3cm} 
   We now turn to stability issues for the above solitary wave, 
   following Angulo \cite{Ang} who solved the problem with $\epsilon=1$ and actually proved using the method in \cite{GSS} that the solitary wave is orbitally stable in $L^2(\R)\times H^1(\R)$ for $|c|<1.$ 
   
   This uses the fact that ${\bf U}_c$ can be viewed as a critical 
   point of the functional $F=\Phi_1-c\Phi_0,$ that is it satisfies
$$(\Phi_1-c\Phi_0)'({\bf U}_c)=0$$
   where 
   $$\Phi_1(v,\eta)=\frac{1}{2}\int_ \R [v^2+\eta^2+v^2\eta+(v_x)^2]dx,$$
   and 
$$\Phi_0(v,\eta)=\int_\R v\eta dx.$$  

Angulo  proved moreover, using the higher order conservation law  
that the Cauchy problem is globally well-posed in $H^{s-1}(\R)\times 
H^s(\R), s\geq 2$ provided the initial  data $(\eta_0,v_0)$ are $L^2(\R)\times H^1(\R)$-close to some translation of the solitary wave. This  is one of the few known results on the global well-posedness of the Cauchy problem for abcd Boussinesq systems.

\begin{remark}
   We are not aware of explicit formulas for  N-solitons of  the good Kaup-Broer-Kupershmidt system.
   
   \end{remark}

\subsection{Remarks on the stationary solutions}

We comment here on the stationary solutions of the good KBK system (here we keep the dependance in $\epsilon$).  Assuming that such a solution vanishes at infinity. In one spatial dimension $(\eta,v)$ should satisfy
\begin{equation}\label{SW-B-Sc-1d}
- v''+ v-\frac{\epsilon^2}{2}v^3 =0,
\end{equation}
$$\eta=-\frac{\epsilon}{2}v^2.$$
and  $v$ can be expressed in terms of  the profile  of the soliton of the focusing cubic nonlinear Schr\"{o}dinger equation,
leading to the explicit form
$$
v(x)=\frac{2}{\epsilon\cosh(\epsilon^{-1/2}x)}.
$$

In the two-dimensional case, one has 
\begin{equation}\label{SW-B-Sc-2d}
\nabla\cdot[-\Delta {\bf v}+{\bf v}-\frac{\epsilon^2}{2}|{\bf v}|^2{\bf v}]=0,
\end{equation}
$$\eta=-\frac{\epsilon}{2}|{\bf v}|^2.$$

Any solution ${\bf v}$ of 
\begin{equation}\label{SW-B-Sc-2d-bis}
-\Delta {\bf v}+{\bf v}-\frac{\epsilon^2}{2}|{\bf v}|^2{\bf v}=0,
\end{equation}
solves  \eqref{SW-B-Sc-2d}.

Equation \eqref{SW-B-Sc-2d-bis} is a particular case of the equation of a bound state solution of a vector nonlinear Schr\"{o}dinger equation. It is proven in \cite{CW} that \eqref{SW-B-Sc-2d-bis} has a solution whose components are constant multiples of the ground state of the corresponding {\it scalar} nonlinear Schr\"{o}dinger equation.

In both cases it would be interesting to investigate the stability of those solutions with respect to the good KBK system.

\section{ The Kaup-Broer-Kupershmidt system by IST techniques}
The connection of \eqref{KBK} with the theory of integrable systems was first noticed in \cite{Kup, Ka} for both the "bad" and "good" cases. Actually, Kaup \cite{Ka},  proved that the KBK system is the compatibility condition for a  pair of linear equations. 

Writing the KBK system as

\begin{equation}\label{Kaupstyle}
\zeta_t+(\zeta v)_x\pm v_{xxx}=0,\quad v_t+vv_x+\zeta_x=0, \quad \zeta=1+\eta,
\end{equation}

 the two aforementioned  linear equations are:

\begin{equation}
    \label{Kauplin}
    \left\lbrace
    \begin{array}{l}
    4\psi_{xx}=\pm[(\lambda-\frac{1}{2}v)^2-\zeta]\psi, \\
    \psi_t=\frac{1}{4}v_x\psi-(\lambda+\frac{1}{2}v)\psi_x.
\end{array}\right.
    \end{equation}
    Kaup  also wrote the formulation of the direct and inverse scattering problems and the formula for the soliton solution.

 Kupershmidt \cite{Kup} showed how the KBK system can be derived in bi-Hamiltonian form. He considered in fact a more general class of systems that he described as the "richest integrable known system known to date":
 
 \begin{equation}
    \label{Kupnew}
    \left\lbrace
    \begin{array}{l}
    u_t+h_x+\beta u_{xx}+uu_x=0, \\
    h_t+\alpha u_{xxx}-\beta h_{xx}=0,
\end{array}\right.
    \end{equation}
    
    where $\alpha$ and $\beta$ are arbitrary real constants. Note that the system is linearly well-posed if and only if $\alpha<-\beta^2.$ The KBK systems correspond to $\beta =0,$ $\alpha >0$ for the "bad" one, $\alpha<0$ for the "good" one.  Kupershmidt derived a corresponding hierarchy by making use of  the theory of non standard integrable systems.
    
    \vspace{0.3cm}
  As already alluded to, the integrable system structure of the KBK systems yields existence of various special solutions, both in the "good" and the "bad" case.

 For instance Sachs \cite{Sachs}   wrote an infinite family of rational solutions of the "bad"  KBK system. Clarkson \cite{Clark} obtains a larger class of rational solutions in terms of the generalized Hermite and generalized Okamoto polynomials for both the "good" and "bad" version of the KBK system.
 
  Furthermore, Ito in \cite{Ito} (see also \cite{Sachs})  exhibited an infinite family of commuting Hamiltonian flows associated to an infinite set of commuting integrals $F_n$ such that 
 
 $$\lbrace F_n,\mathcal H\rbrace=0$$
 
 where $\lbrace, \rbrace$ is the Poisson bracket
 and $\mathcal H$ the Hamiltonian of the KBK system.

  This leads to  a hierarchy of KBK Boussinesq systems.
 
 \vspace{0.3cm}
 
  The second one in the "bad" case  is
 
 \begin{equation}
    \label{KBK2nd}
    \left\lbrace
    \begin{array}{l}
    v_t+\frac{1}{4}(v^3+6v\zeta+4v_{xx})_x=0, \\
    \zeta_t+\frac{1}{4}(3v^2\zeta+3\zeta^2+3v^2_x+6vv_{xx}+4\zeta_{xx})_x=0,
\end{array}\right.
 \end{equation}
 
 where again $\zeta=1+\eta.$
 
 Note that, contrary to the "bad" KBK system, the linear part in \eqref{KBK2nd} is decoupled and well-posed.
 
 The  Hamiltonian corresponding  to \eqref{KBK2nd} is
 
 \begin{equation}\label{Ham2nd}
 H_2(v,\zeta)=\frac{1}{4}\int_\R\lbrace v^3\zeta+3v\zeta^2+2v_{xx}\zeta+2v\zeta_{xx}+\frac{3}{2}v^2v_{xx}\rbrace dx.
 \end{equation}
 
 \begin{remark}
 The system \eqref{KBK2nd} is different from another integrable system introduced by Kupershmidt \cite{Kup2} as the  {\it s-mKdV system} and which writes, in the notations of \cite{ZTL} where one also finds a bi-Hamiltonian formulation of the equation:
 
  \begin{equation}\label{s-mKdV}
    \left\lbrace
    \begin{array}{l}
    v_t-v_{xxx}+6v^2v_x-3(\eta_{xx}\eta)_x-6(\eta_x\eta v)_x=0, \\
    \eta_t-4\eta_{xxx}+6(v^2-v_x)+3(2v_xv-v_{xx})\eta=0.
\end{array}\right.
 \end{equation}
 \end{remark}
 
 Actually \eqref{s-mKdV} is obtained after  a Miura type transform 
 $$u=-v_x-v^2+\eta_x\eta,\quad \xi=\eta_x+\eta v$$
 
 from the so-called super-KdV equation (s-KdV), see \cite{Kup2, ZTL}
 
 \begin{equation}\label{s-KdV}
 \left\lbrace
    \begin{array}{l}
    u_t-u_{xxx}-6uu_x+12\xi\xi_{xx}=0, \\
    \xi_t-4\xi_{xxx}-6u\xi_x-3\xi u_x=0.
\end{array}\right.
\end{equation}

  We refer to \cite{Pil, Bar, PM} for ill-posedness and well-posedness results on the Cauchy problem for a variant of  \eqref{s-KdV}. On the other hand we are not aware of mathematical results on the Cauchy problem for \eqref{s-mKdV}.

 \vspace{0.5cm}
 The next system in the "bad" KBK  hierarchy  is in the $(v,\zeta)$ variables
 
  \begin{equation}
    \label{KBK3rdd}
    \left\lbrace
    \begin{array}{l}
    v_t+\frac{1}{8}(v^4+12 v^2\zeta+6\zeta^2+6v_x^2+16vv_{xx}+8\zeta_{xx})_x=0, \\
    \zeta_t+\frac{1}{8}(4v^3\zeta+12v\zeta^2+16v\zeta_{xx}+20v_x\zeta_x+20v_{xx}\zeta+8v_{xxxx}+12vv_x^2+12v^2v_{xx})_x=0,
\end{array}\right.
 \end{equation}
 which is linearly ill-posed.

 \vspace{0.3cm}
 More generally the higher order flows are defined by the equations
 
 \begin{equation}\label{Hierarchy}
 \begin{pmatrix}v
\\\zeta\end{pmatrix}_t+J\nabla H_m=0,
 \end{equation}
 
 where 
 
 $$J=\begin{pmatrix}0&\partial_x\\\partial_x&0\end{pmatrix},$$
 
 and where the Hamiltonian $H_m$ are defined inductively by
 
 $$H_2=\frac{1}{2}\int_\R(v^2\zeta+\zeta^2-v_x^2)dx,$$
 
 $$J\nabla H_m=K\nabla H_{m-1},$$
 where K is the skew-symmetric operator

  $$K=\begin{pmatrix}\partial_x&\frac{1}{2}\partial_x v\\
 \frac{1}{2}v\partial_x&\partial_x^3+\zeta\partial_x+\frac{1}{2}\zeta_x\end{pmatrix}.$$
 
 \vspace{0.3cm} 
 The same strategy can be applied to the "good" KBK system, starting now from the Hamiltonian
 
 $$H_2=\frac{1}{2}\int_\R(v^2\zeta+\zeta^2+v_x^2)dx,$$
 
 and the new K defined  as
 
 $$K=\begin{pmatrix}\partial_x&\frac{1}{2}\partial_x v\\
 \frac{1}{2}v\partial_x&-\partial_x^3+\zeta\partial_x+\frac{1}{2}\zeta_x\end{pmatrix}.$$
 
The second equation in the hierarchy is thus:

\begin{equation}
    \label{goodKBK2nd}
    \left\lbrace
    \begin{array}{l}
    v_t+\frac{1}{4}(v^3+6v\zeta-4v_{xx})_x=0, \\
     \zeta_t+\frac{1}{4}(3v^2\zeta+3\zeta^2-3v_x^2-6vv_{xx}-4\zeta_{xx})_x=0,
\end{array}\right.
 \end{equation}
which again is linearly well-posed and correspond to the Hamiltonian

 \begin{equation}\label{Ham2ndgood}
 H_2(v,\zeta)=\frac{1}{4}\int_\R\lbrace v^3\zeta+3v\zeta^2-2v_{xx}\zeta-2v\zeta_{xx}-\frac{3}{2}v^2v_{xx}\rbrace dx.
 \end{equation}
% Assuming that $v$ and its derivatives tend to $0$ as $|x|\to \infty,$ one can write
 
% \begin{equation}\label{Ham2ndgood-bis}
% H_2(v,\zeta)=\frac{1}{4}\int_\R\lbrace v^3\zeta+3v\zeta^2+4v_x^2\zeta_x^2+3vv_x^2\rbrace dx.
% \end{equation}
 
\begin{remark}
In terms of the variables $ (v, \eta),$\eqref{goodKBK2nd} writes 
\begin{equation}
    \label{goodKBK2nd(2)}
    \left\lbrace
    \begin{array}{l}
    v_t+\frac{1}{4}(v^3+6v+6v\eta-4v_{xx})_x=0, \\
     \zeta_t+\frac{1}{4}(3v^2+3v^2\eta+6\eta+3\eta^2-3v_x^2-6vv_{xx}-4\eta_{xx})_x=0,
\end{array}\right.
 \end{equation}

\end{remark}

 \begin{remark}
  Equations \eqref{KBK2nd},  \eqref{goodKBK2nd} and \eqref{goodKBK2nd(2)} are systems of KdV type involving second and third order nonlinear terms making the study of the local Cauchy problem delicate. We refer to \cite{CKS, LPS} for a study of the Cauchy problem for scalar  dispersive equations involving a nonlinear dispersive third order term. 
 \end{remark}

 The next system in the "good" KBK hierarchy is 
 
 \begin{equation}
    \label{goodKBK3rd}
    \left\lbrace
    \begin{array}{l}
    v_t+\frac{1}{8}(v^4+12v^2\zeta +6\zeta^2-6v_x^2-16vv_{xx}-8\zeta_{xx})_x=0, \\
     \zeta_t+\frac{1}{8}(4v^3\zeta+12v\zeta^2-16v\zeta_{xx}-20v_x\zeta_x-20v_{xx}\zeta+8v_{xxxx}-12vv_x^2-12v^2v_{xx})_x=0,
\end{array}\right.
 \end{equation}
 
 which is linearly well-posed and the $(v,\zeta)$ version of \eqref{KBKJaime}.
 
 \vspace{0.3cm}
It corresponds to the Hamiltonian
 \begin{equation}\label{newHam}
 H(v,\zeta)=\frac{1}{8}\int_\R (v^4\zeta+6v^2\zeta^2-8v^2\zeta_{xx}+4v^2_{xx}+10\zeta v_x^2+4\zeta_x^2-2v^3v_{xx})dx.
 \end{equation}

 %and it is the $(v,\zeta)$ version of \eqref{KBKJaime}.
 
 %\begin{remark}
 %We are not aware of the physical relevance of the higher order KBK systems.
 %\end{remark}
 
\vspace{0.3cm}

We now briefly describe some recent progress on the Kaup-Broer-Kupershmidt  system. 
   % It is recalled in \cite{NZ}  that  Broer \cite{Broer} and Kaup \cite{Ka} derived the system in dimensional form.
  %Kaup has shown that \eqref{KBK} can be solved by the Inverse 
  %Scattering Method, Kupershmidt proving the complete integrability. Nevertheless, considering that this Boussinesq system is ill-posed, this property seems irrelevant to solve the Cauchy problem, but it allows the construction of one and N-solitons, periodic solitons, undular bores,  see for instance \cite{EGP, Ka, NZ}. On the other hand the case with strong surface tension, that is $\tau/(g\rho)>h_0^2/3$ corresponds to $b=d=c=0, a<0$ in \eqref{abcd} and thus to the good KBK equation. It also have an integrable structure, for instance Kaup \cite{Ka} found the Lax pair,  Matveev and Yavor \cite{MY} computed  finite gap solutions (also in the case of small surface tension).
   
   The paper \cite{NAF}  focusses on the periodic "good" KBK system and provides an explicit form of the periodic traveling waves in terms of the Weierstrass elliptic function $\wp:$
  
  $$v(\xi)=\frac{a_{11}\wp(\xi+\phi)+a_{12}}{a_{21}\wp(\xi+\phi)+a_{22}},\quad \eta(\xi)=cv(\xi)-\frac{1}{2}v(\xi)^2;\;\;\ \xi=(x-ct).$$

  The authors  also exhibit a matrix  Lax pair   that is used it to prove the existence of an infinite number of complex conservation laws:
    \begin{theorem}-\cite{NAF}

  $$(\rho_n)_t+(j_n)_x=0,$$
  
  where the conserved densities $\rho_n$ are determined by the recursion:
  
  $$\rho_1=\frac{1}{2}\eta+\frac{i}{2}v_x,\quad \rho_2=iv\rho1-2\rho_{1x}+\frac{i}{2}v,$$
  
  $$\rho_{n+1}=iv\rho_n-\rho_{n-1}-2\rho_{nx}-2\sum_{k=1}^{n-1}\rho_k\rho_{n-k},\quad n>1,$$
  
  and the conserved currents are determined by the conserved densities:
  
  $$j_1=\frac{1}{4}v+\frac{1}{2}v\rho_1-\frac{i}{2}\rho-2,$$
  
  $$j_n=\frac{1}{2}v\rho_n-\frac{i}{2}\rho_{n+1}+\frac{i}{2}\rho_{n-1},\quad n>1.$$
  
  Each complex conservation law gives rise to two real conservation laws that are polynomials in $\eta,v$ and their higher order derivatives with respect to x.
  
  \end{theorem}
  
  The Lax pair is also used to construct the forward scattering transform for periodic solutions and to find exact formulas for finite gap solutions.

\vspace{0.3cm}
Two-dimensional integrable generalizations of the 
Kaup-Broer-Kupershmidt system are given in \cite{RP, NZ, Nab}. None of them seem physically relevant, however.

On the other hand for capillary waves such that $\tau/(g\rho)>h_0^2/3,$ \eqref{BKK} is analogous to the good BSK system.

\begin{remark}
The "bad" Kaup-Broer-Kuperschmidt system is studied from the Inverse Scattering point of view in  \cite{ZW} under the form
\begin{equation}
    \label{BZW}
    \left\lbrace
    \begin{array}{l}
    v_t=\frac{1}{2}(v^2+2w-v_x)_x, \\
    w_t=(vw+\frac{1}{2}w_x)_x,
\end{array}\right.
    \end{equation}
which is an alternative form of \eqref{new-BKK}. Using this version of the system, \eqref{BZW} is the compatibility condition of the following Lax pair
\begin{equation}\label{BLax}
\psi_x=U\psi,\quad \psi_t=V\psi,
\end{equation}
with 
\begin{equation}
    \label{BLax2}
    \left\lbrace
    \begin{array}{l}
    U=-ik\sigma_3+Q,\quad Q=\frac{v}{2}\sigma_3+\sigma_+-w\sigma_-, \\
    V=k^2\sigma_3+\widehat{Q},\quad \widehat{Q}=-\frac{1}{4}(v_x-v^2)\sigma_3+(ik+\frac{v}{2})\sigma_+-\frac{1}{2}(w_x+wv)\sigma_-
\end{array}\right.
    \end{equation}
where $k$ is a spectral parameter, $\sigma_j, j=1,2,3$ are the classical Pauli matrices and $\sigma_+=\frac{1}{2}(\sigma_1+i\sigma_2),\; \sigma_-=\frac{1}{2}(\sigma_1-i\sigma_2).$

The Inverse Scattering mechanism described in \cite{ZW} allows the construction of special solutions such as kinks, see \cite{ZW} for details.
\end{remark}

\begin{remark}
There is so far no rigorous result on the complete resolution of the 
Cauchy problem for the well-posed Kaup-Broer-Kupershmidt system 
($\alpha=-1$) under an appropriate functional setting for the initial 
data and the solution by using inverse scattering techniques, and 
giving possibly insight into the qualitative behavior of the solutions.This would include a rigorous theory for the direct and inverse scattering problem. We refer to \cite{NAF} for progress in this direction.
\end{remark}

\vspace{0.3cm}

 \section{Numerical approach to the KBK system}
In this section, we will detail the numerical approach used to solve the 
`good' KBK system in the numerical experiments. 

We recall the conserved energy for this equation is
\begin{equation}
	E = \frac{1}{2}\int_{\mathbb{R}}^{}  
	(\eta^{2}+(1+\eta)v^{2}+v_{x}^{2})dx
	\label{E}.
\end{equation}
Solitons of the KBK system (with location $x_{0}$ of the maximum for $t=0$)
can be written for the velocity 
$C\in\mathbb{R}$, $|C|<1$ with (\ref{SWexpl}) in the form
\begin{align}
	v& =\frac{2(1-C^2)}{\cosh(\sqrt{1-C^2}(x-Ct-x_{0}))-C},
	\nonumber\\
	\eta & = Cv-\frac{1}{2}v^2 
	\label{soliton}.
\end{align}

To numerically solve the system (\ref{KBK}) we essentially use the diagonalisation 
approach for the linear part of (\ref{bis}), (\ref{ter}). We consider the 
KBK system in Fourier space,
\begin{align}
	\hat{\eta}_{t} & = - ik(1+k^{2})\hat{v}-ik\widehat{\eta v},
	\nonumber\\
	\hat{v}_{t} & =-ik \hat{\eta}-\frac{ik}{2}\widehat{v^{2}}
	\label{KBKfourier}.
\end{align}
Here we use  the standard definition for a  Fourier transform for 
integrable
functions $u(x)$ 
denoted by $\hat{u}(k)$ with dual variable $k$ 
and its inverse (in the sense of tempered distributions),
\begin{align}
	\hat{u}(k)& = \int_{\mathbb{R}}^{}u(x)e^{-ikx}dx,\quad k\in 
	\mathbb{R},
	\nonumber\\
	u(x) & =\frac{1}{2\pi}\int_{\mathbb{R}}^{}\hat{u}(k)e^{ikx} 
	dk,\quad x\in\mathbb{R}.
	\label{fourierdef}
\end{align}

Introducing 
\begin{equation}
	\hat{u}_{\pm}= \hat{v}\pm\frac{\hat{\eta}}{\sqrt{1+k^{2}}}
	\label{hatu},
\end{equation}
we can write the system (\ref{KBKfourier}) in the form 
\begin{equation}
	(u_{\pm})_{t}=\mp u_{\pm}-ik\left(\frac{1}{2}\widehat{v^{2}} 
	\pm\frac{\widehat{\eta v}}{\sqrt{1+k^{2}}}\right)
	\label{upm}.
\end{equation}
Obviously the dispersion relation is as in (\ref{new-BKK2}), i.e., 
for large $|k|$ as in the Schr\"odinger equation.

For the numerical treatment, the Fourier transform in (\ref{upm}) 
will be approximated in standard 
way by the discrete Fourier transform (DFT) which can be conveniently 
computed by the \emph{Fast Fourier Transform} (FFT). This is a 
\emph{spectral method} which means that the numerical error in 
approximating smooth periodic functions with $N$ modes in the DFT 
leads to a numerical error exponentially decreasing with $N$, see the 
discussion in \cite{trefethen} and references therein. Thus we will 
in the following always work on a torus of period $2\pi L$ with 
$L>0$, i.e., we will consider values of $x\in L[-\pi,\pi]$, where we 
will apply $N$ DFT modes. In an abuse of notation, we will denote the 
discrete Fourier transform with the same symbol as the standard 
Fourier transform. Since the numerical error in approximating a 
function with a DFT is of the order of the highest DFT coefficients, 
we always use values of $L$ and $N$ such that the DFT coefficients 
decrease to machine precision which is here of the order of 
$10^{-16}$. 

The resulting $2N$ dimensional system of equations (\ref{upm}) are of the form 
\begin{equation}
	u_{t}=\mathcal{L}u + \mathcal{N}(u)
	\label{LN},
\end{equation}
where $\mathcal{L}$ is a linear diagonal operator, here proportional 
to $\pm ik\sqrt{1+k^{2}}$, whereas $\mathcal{N}(u)$ is a nonlinear 
term in the $u_{\pm}$. Due to  $|\mathcal{L}|$ being rapidly 
increasing with $|k|$, the system is \emph{stiff} which means that 
explicit time integration schemes are not efficient for such systems 
due to stability conditions, 
see for instance the discussion in \cite{HO} and references therein. 
However, there are efficient time integration schemes to address 
systems of the form (\ref{LN}) with a stiff diagonal linear term, see 
\cite{HO} for so-called exponential time differencing schemes (ETD). 
The idea of ETD schemes is to use equidistant time 
steps $h$ and to integrate equation (\ref{LN})  between 
the time steps $t_{n}$ and $t_{n+1}$, $n=1,2,\ldots$ with an 
exponential integrator with respect to 
$t$. We get 
$$
u(t_{n+1})=e^{\mathcal{L}h}u(t_{n})+\int_{0}^{h} e^{\mathcal{L}(h-\tau)}
    \mathcal{N}(u(t_{n}+\tau),t_{n}+\tau)d\tau.
$$
The integral will be computed in an approximate way for which 
different schemes exist. In \cite{KR} we have compared for dispersive 
PDEs various   
Runge-Kutta schemes of classical order 4 which all showed similar 
performance. Therefore we apply in the following the one by   Cox-Matthews 
\cite{CM}. Since the method is explicit (only information of the 
functions at the previous time step $t_{n}$ is needed), it is not important that 
the nonlinear part is not diagonal in (\ref{upm}) as the linear part. 
As discussed in \cite{etna,KR}, the exactly conserved energy can be 
used to control the numerical error in the time integration. Whereas 
the energy (\ref{E}) is exactly conserved for the KBK system, the 
numerically computed energy will depend on time due to unavoidable 
numerical errors. The relative energy
\begin{equation}
	\Delta :=|E(t)/E(0)-1|
	\label{rel}
\end{equation}
typically overestimates the numerical error by 1-2 orders of 
magnitude. We will always aim at a $\Delta$ considerably smaller than 
$10^{-3}$. 

\begin{figure}[htb!]
\includegraphics[width=0.49\textwidth]{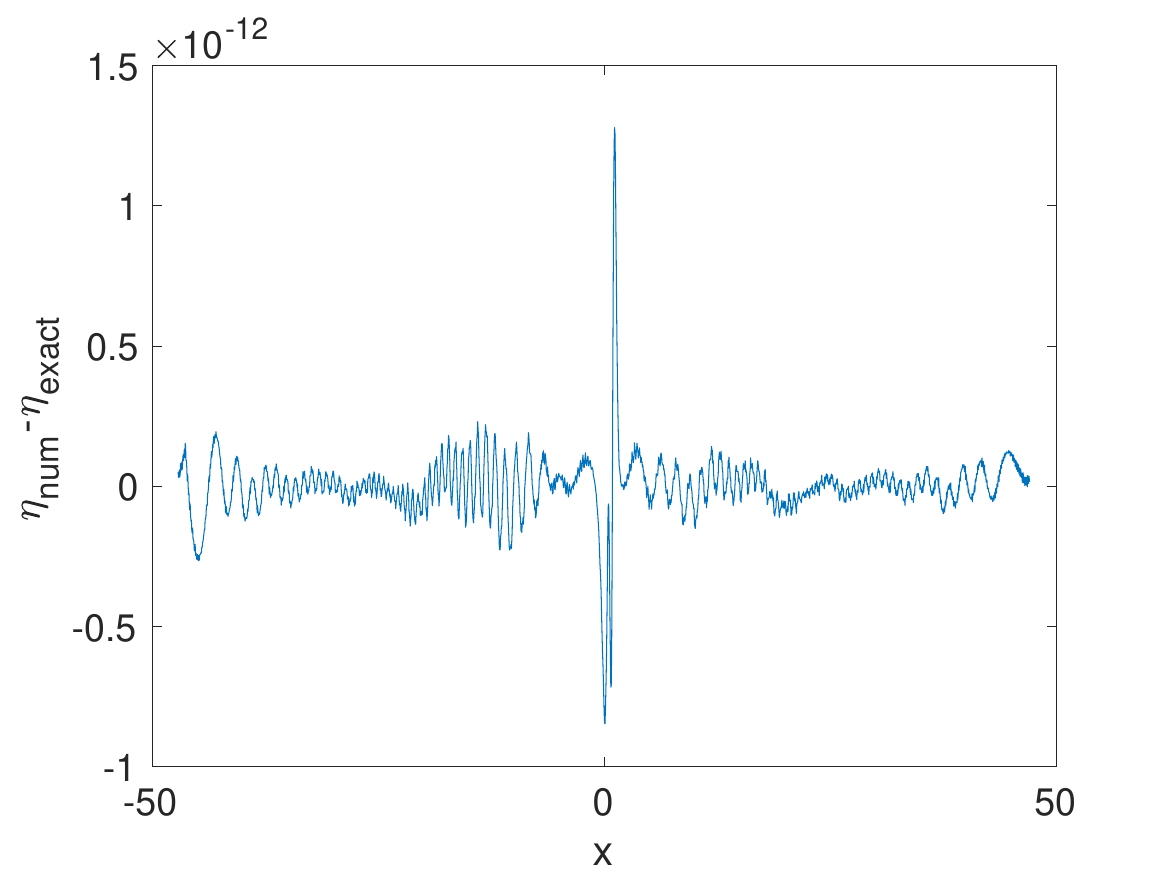}
\includegraphics[width=0.49\textwidth]{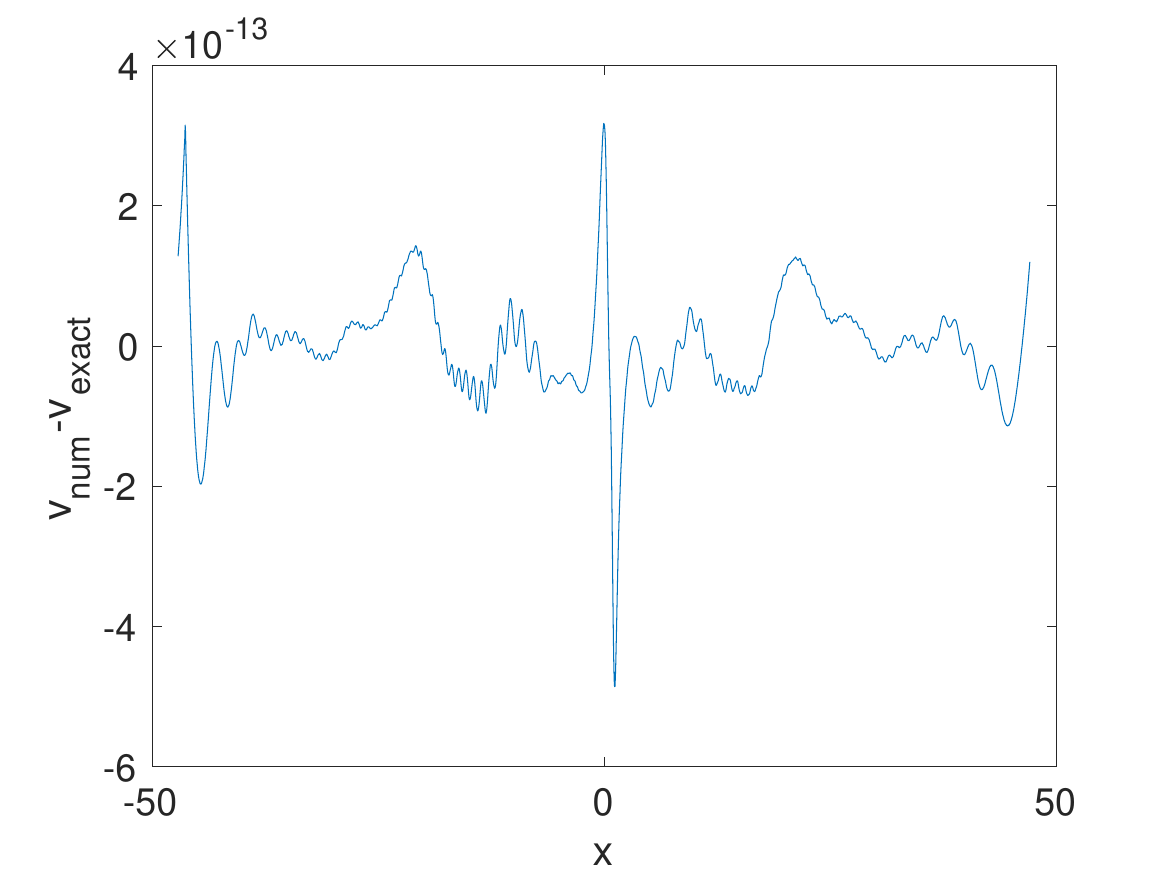}
\caption{Difference between numerical and exact solution for soliton 
initial data (\ref{soliton}) with $C=0.8$ for $t=1$, on the left $\eta$, on 
the right $v$. }
\label{soltest}
\end{figure}

To test the code, we consider the soliton (\ref{soliton}) for $C=0.8$ 
as initial data. We use $N=2^{11}$ DFT modes for $x\in15[-\pi,\pi]$ 
with $N_{t}=4000$ time steps for $t\leq 1$. The DFT coefficients 
decrease to machine precision during the whole computation in this 
case, and the relative energy is conserved to better than $10^{-12}$. 
In Fig.~\ref{soltest} we show the difference between numerical and 
exact solution at the final time. It can be seen that it is of the 
order of $10^{-12}$ which shows that the code can propagate the 
soliton with essentially machine precision. 

\section{Perturbed solitons}
Angulo \cite{Ang} showed that the solitons of the KBK system are 
stable. In this section we illustrate this by considering 
perturbations of the solitons. 

We study perturbations of the form 
\begin{equation}
	v(x,0) = \lambda v_{C}(x),\quad \eta(x,0) = \mu \eta_{C}(x),
	\label{pert}
\end{equation}
where $v_{C}$, $\eta_{C}$ are the solitons (\ref{soliton}) for a given 
real velocity $C$ with $|C|<1$, and where $\lambda$, $\mu$ are real 
numbers in the vicinity of 1. We use $N=2^{12}$ DFT modes for $x\in 
30[-\pi,\pi]$ and $N_{t}=4000$ time steps for $t\in[0,5]$. Note that 
the maximum of $v$ in (\ref{soliton}) is given by $1+C$ and thus 
strictly smaller than 2. This means that the value of $\lambda$ in 
(\ref{pert}) has to be chosen such that this condition is satisfied. 
To study numerically perturbations and to see an effect in finite 
time, one has to consider values of $\lambda$ and or $\mu$ that have 
a finite difference to 1. This implies that the resulting solution 
will be close to the original soliton, but not of identical velocity. 
We fit the soliton in 
the following way: in the numerical solution for $v$ at the final 
time $t=5$, we identify the location $x_{0}$ and the value $v_{0}$ of the 
maximum. Since the maximum value for $v$ in (\ref{soliton}) 
is given by $2(C+1)$, we can get the value of $C$ from this maximum. 

\begin{figure}[htb!]
 \includegraphics[width=0.49\textwidth]{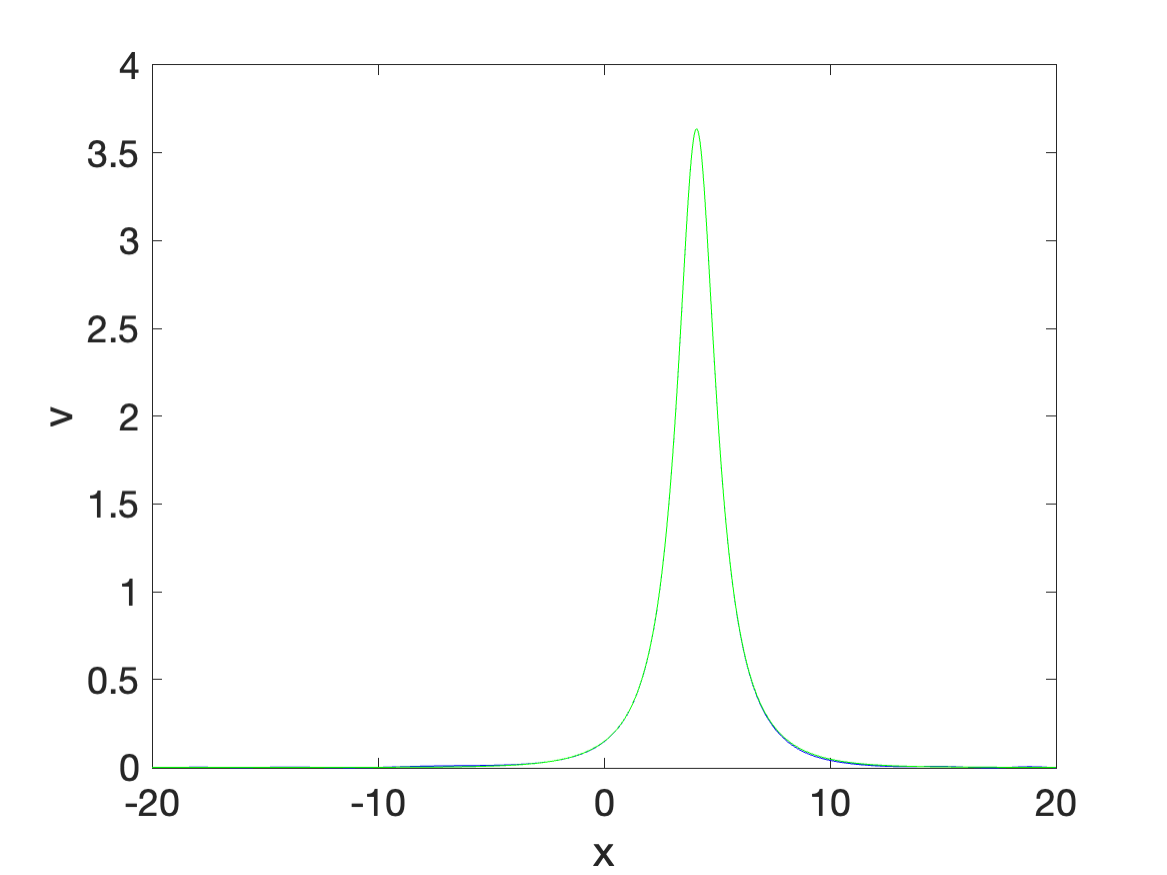}
 \includegraphics[width=0.49\textwidth]{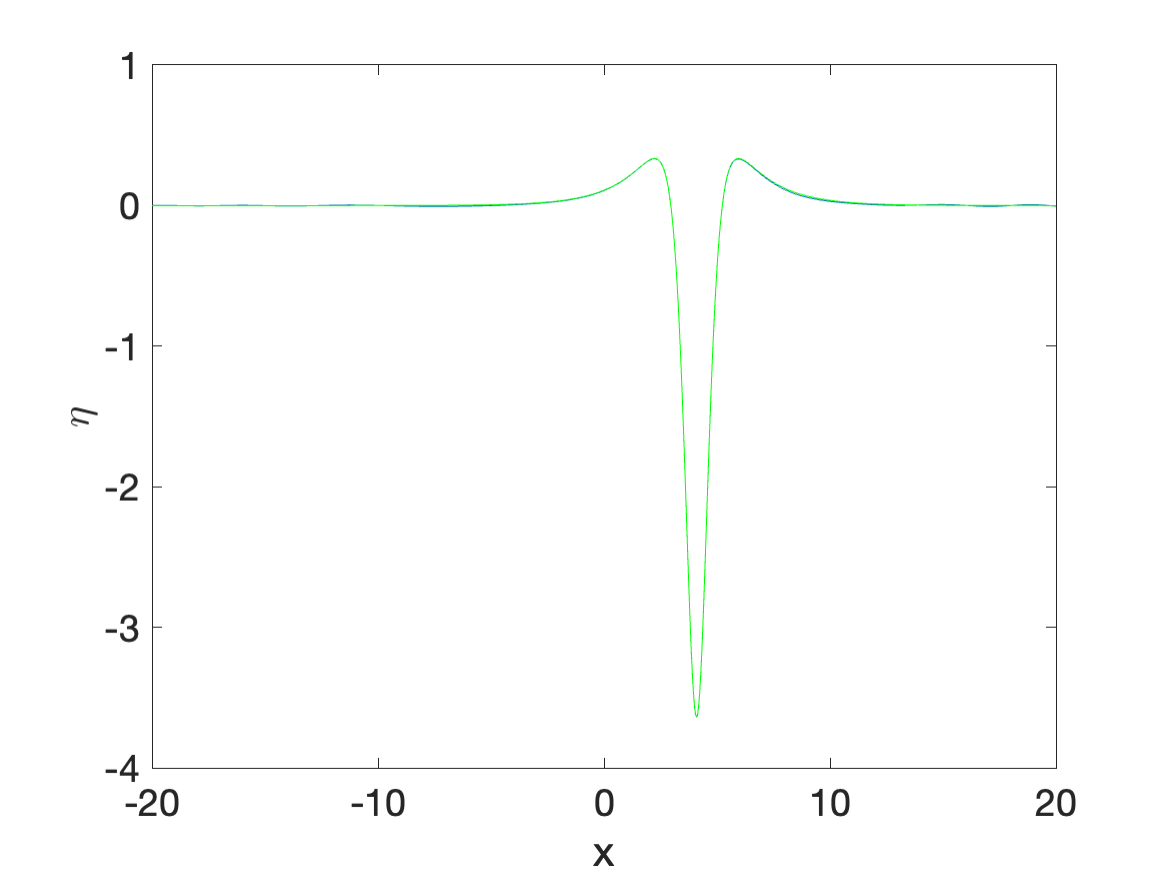}
 \caption{The solution to the KBK system for perturbed  soliton 
 initial data of the form (\ref{pert}) with $\lambda=1.01$ and $\mu=1$
 at the final time $t=5$ in blue and a 
 fitted soliton in green, on the left $v$, on the right $\eta$.  }
 \label{KBKsolc08_101t5}
\end{figure}

We show the solution to the KBK system for $C=0.8$, $\lambda=1.01$ 
and $\mu=1$ at the final time together with the fitted soliton in green (the fitted 
velocity is 0.819) in Fig.~\ref{KBKsolc08_101t5}. It can be seen that 
the solution is very close to the fitted soliton, but that there is 
also some small radiation. The soliton is thus as expected stable. 

\begin{figure}[htb!]
 \includegraphics[width=0.49\textwidth]{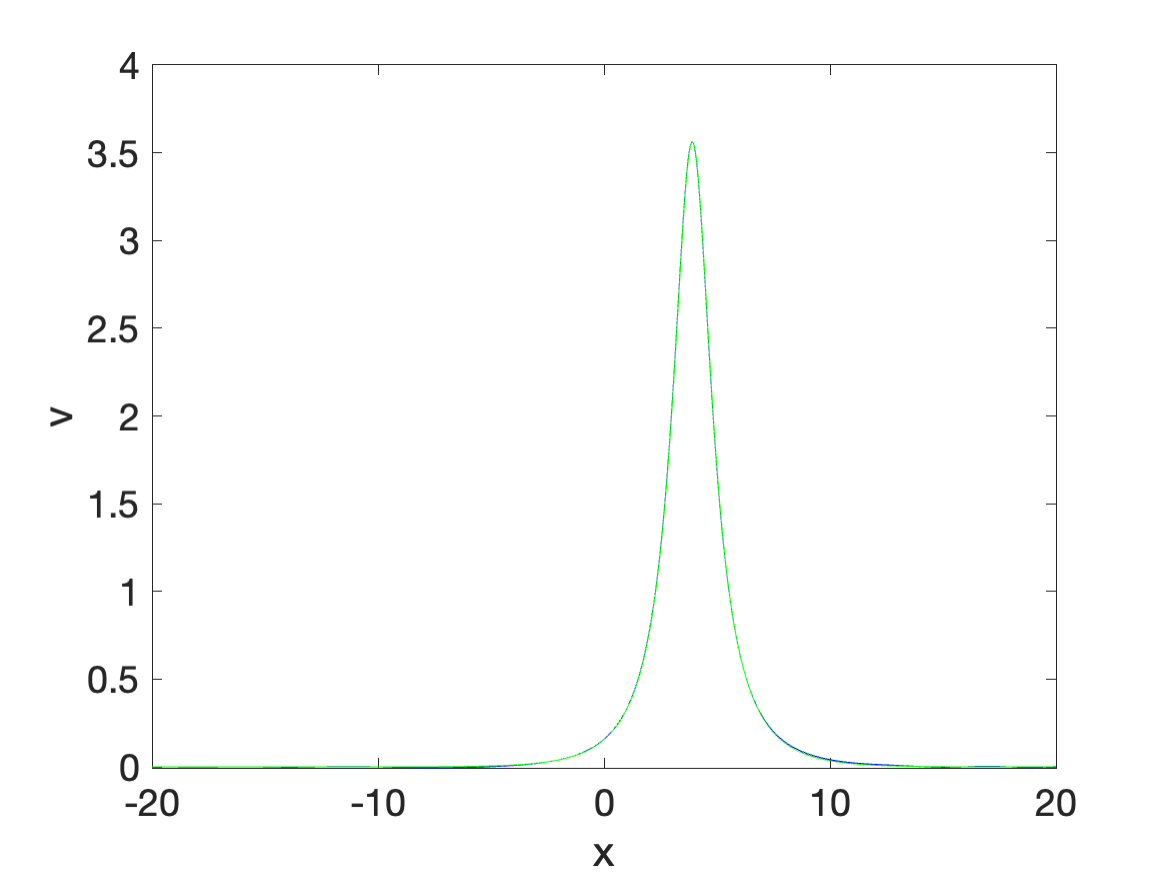}
 \includegraphics[width=0.49\textwidth]{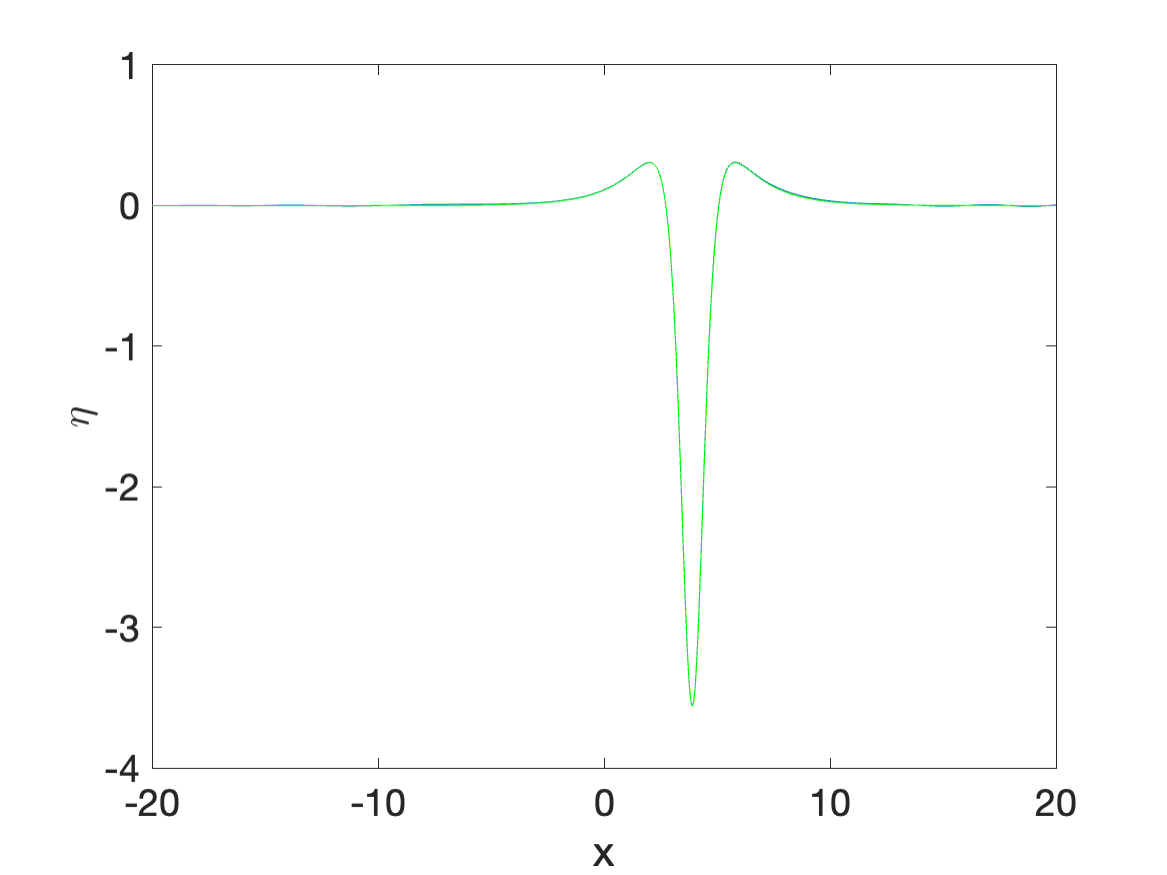}
 \caption{The solution to the KBK system for perturbed  soliton 
 initial data of the form (\ref{pert}) with $\lambda=0.99$ and $\mu=1$
 at the final time $t=5$ in blue and a 
 fitted soliton in green, on the left $v$, on the right $\eta$.  }
 \label{KBKsolc08_099t5}
\end{figure}

The situation is very similar for initial data of the form 
(\ref{pert}) with $\lambda=0.99$ and $\mu=1$ as can be seen in 
Fig.~\ref{KBKsolc08_099t5}. The fitted soliton has velocity $0.7811$. 
Again the final state is a soliton plus radiation.  

\begin{figure}[htb!]
 \includegraphics[width=0.49\textwidth]{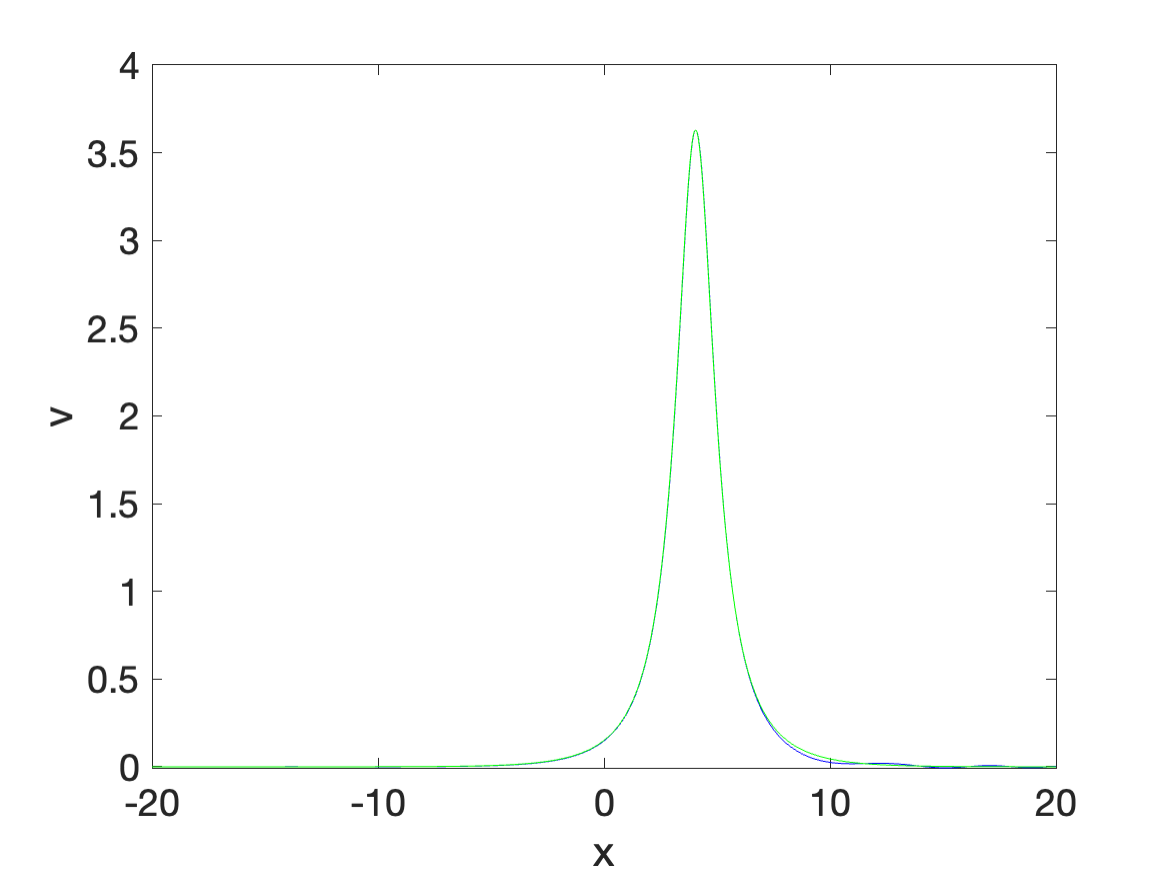} 
 \includegraphics[width=0.49\textwidth]{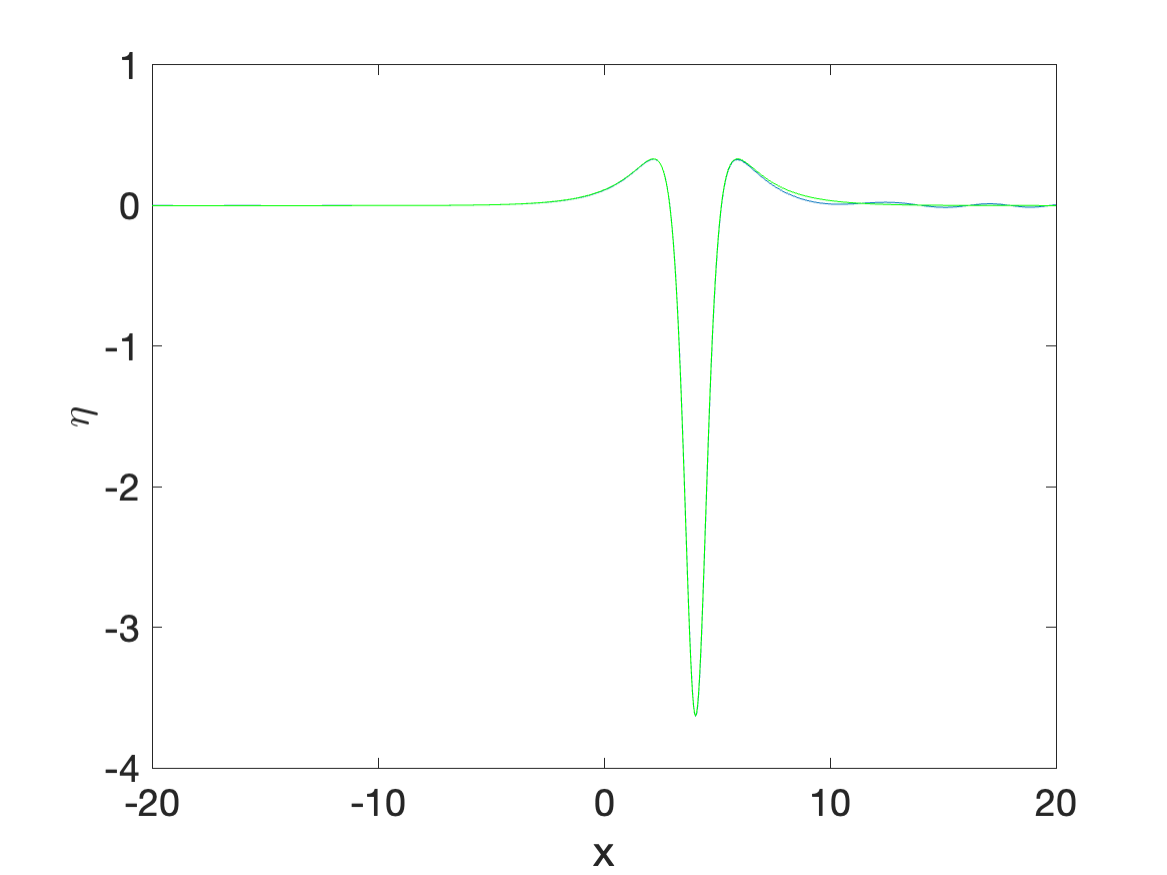}\\
  \includegraphics[width=0.49\textwidth]{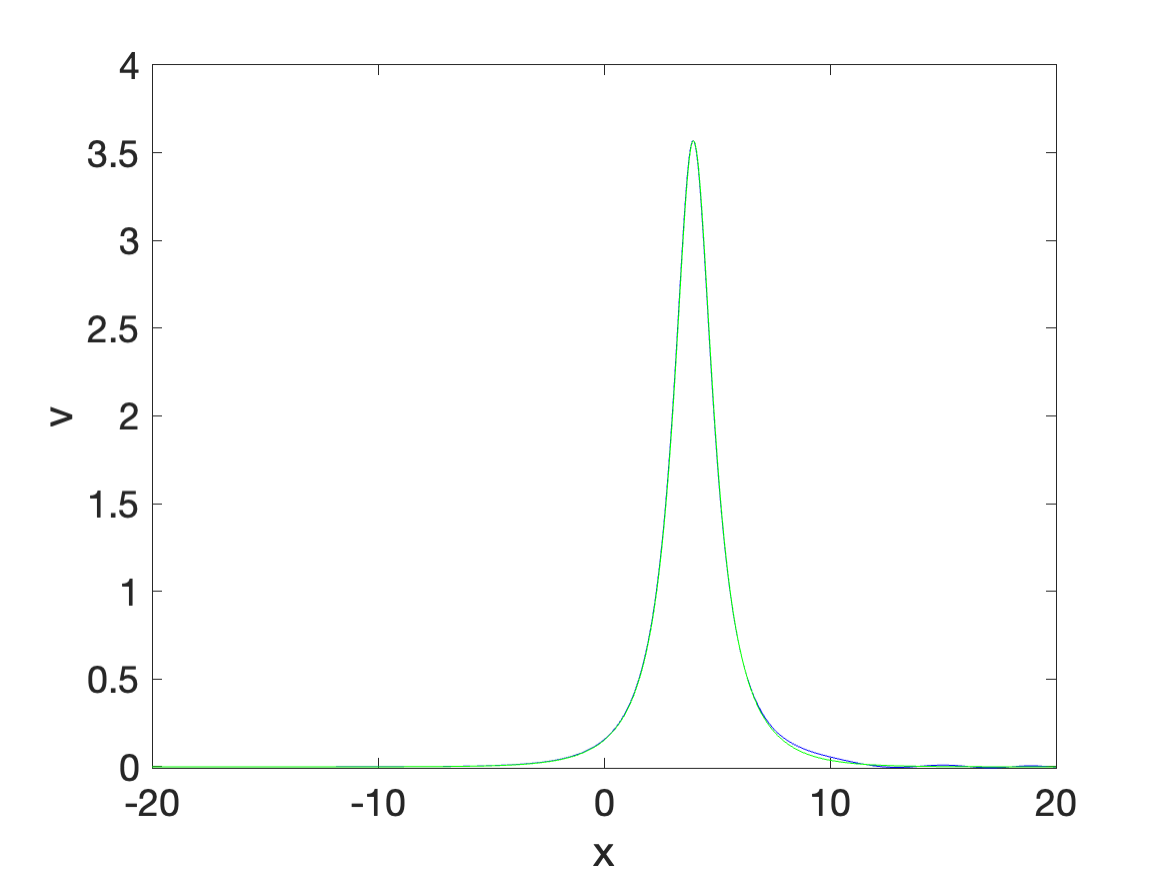}
 \includegraphics[width=0.49\textwidth]{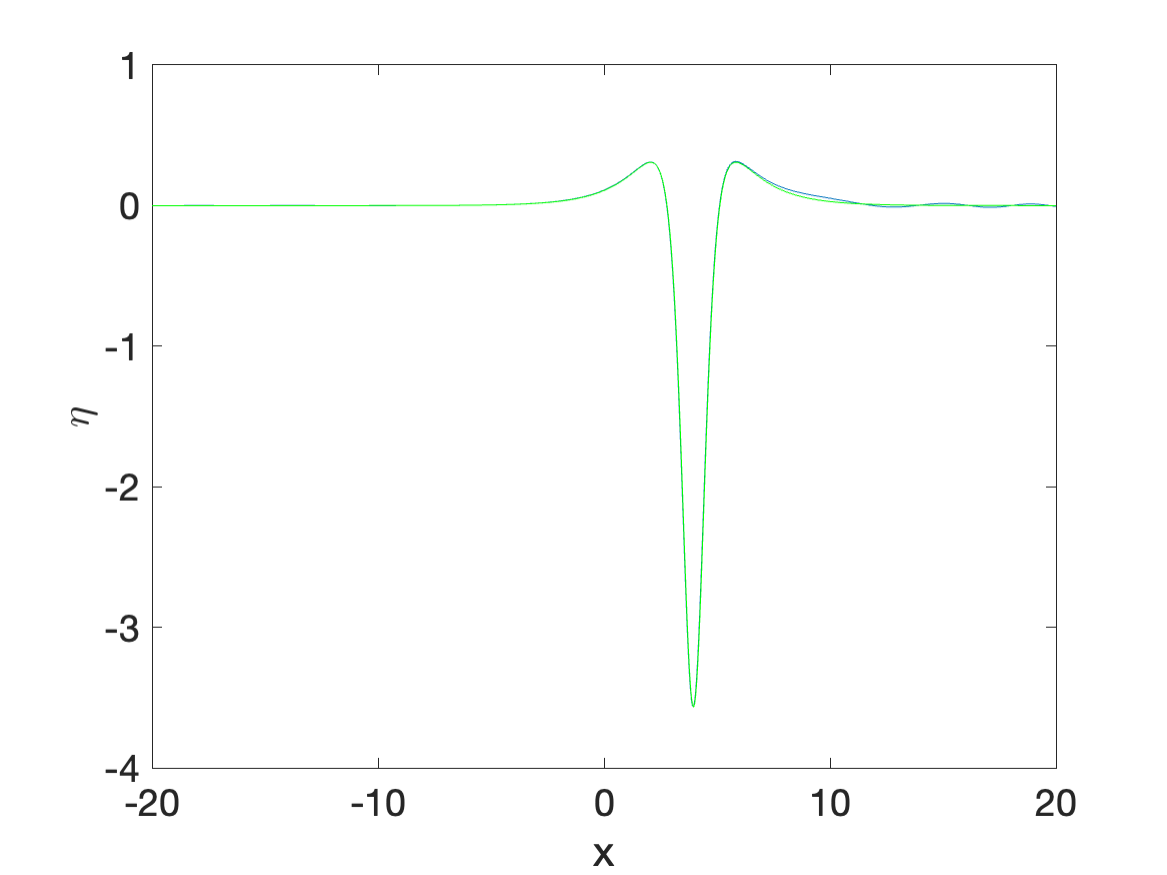}
\caption{The solution to the KBK system for perturbed  soliton 
 initial data of the form (\ref{pert}) with $\lambda=1$ and 
 $\mu=1.01$ in the upper row and $\mu=0.99$ in the lower row 
 fitted soliton in green, on the left $v$, on the right $\eta$.  }
 \label{KBKsolc08_mut5}
\end{figure}

The same stability aspects are observed if perturbations of the form 
(\ref{pert}) with $\lambda=1$ and $\mu\sim1$ are studied as shown in 
Fig.~\ref{KBKsolc08_mut5}. The fitted values of the velocity are 
$0.8155$ for $\mu=1.01$ and $0.7846$ for $\mu=0.99$. In all cases 
the final state of the solution is a soliton plus radiation. 

\begin{figure}[htb!]
 \includegraphics[width=0.49\textwidth]{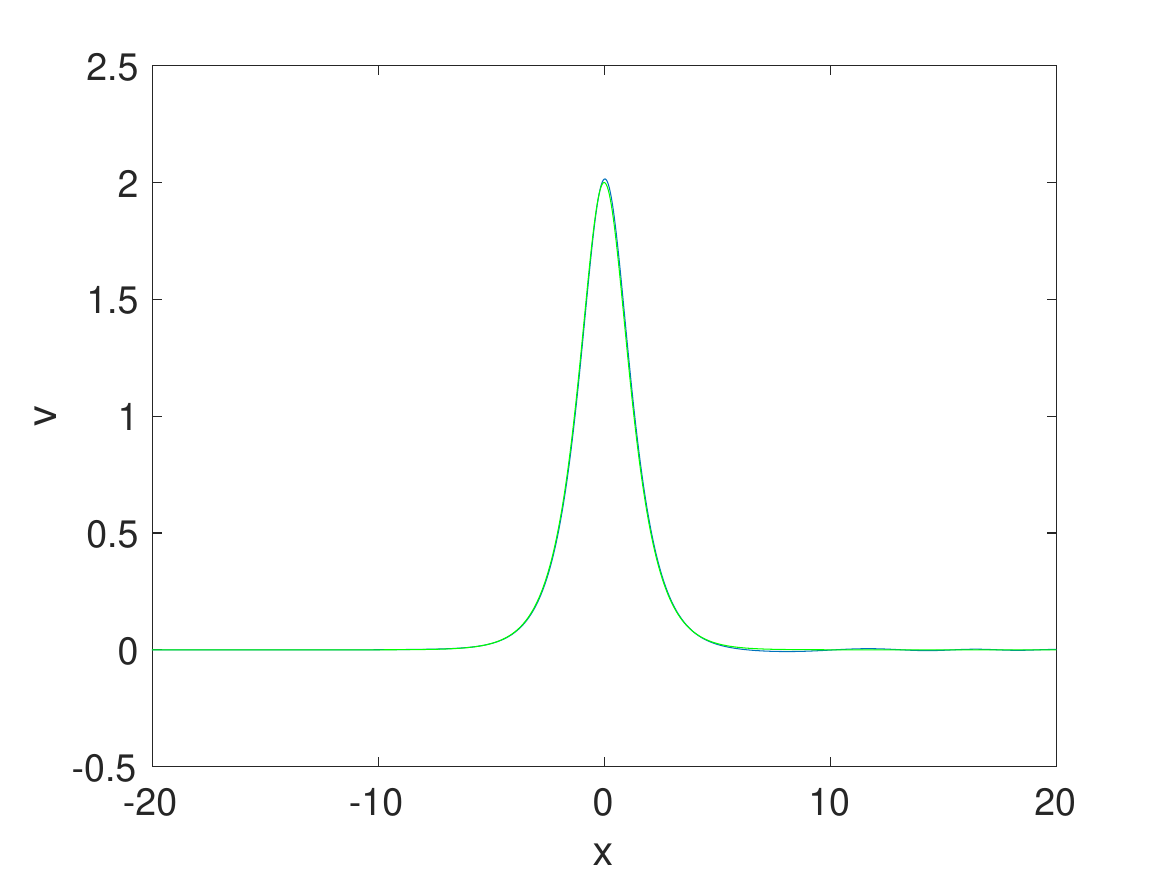}
 \includegraphics[width=0.49\textwidth]{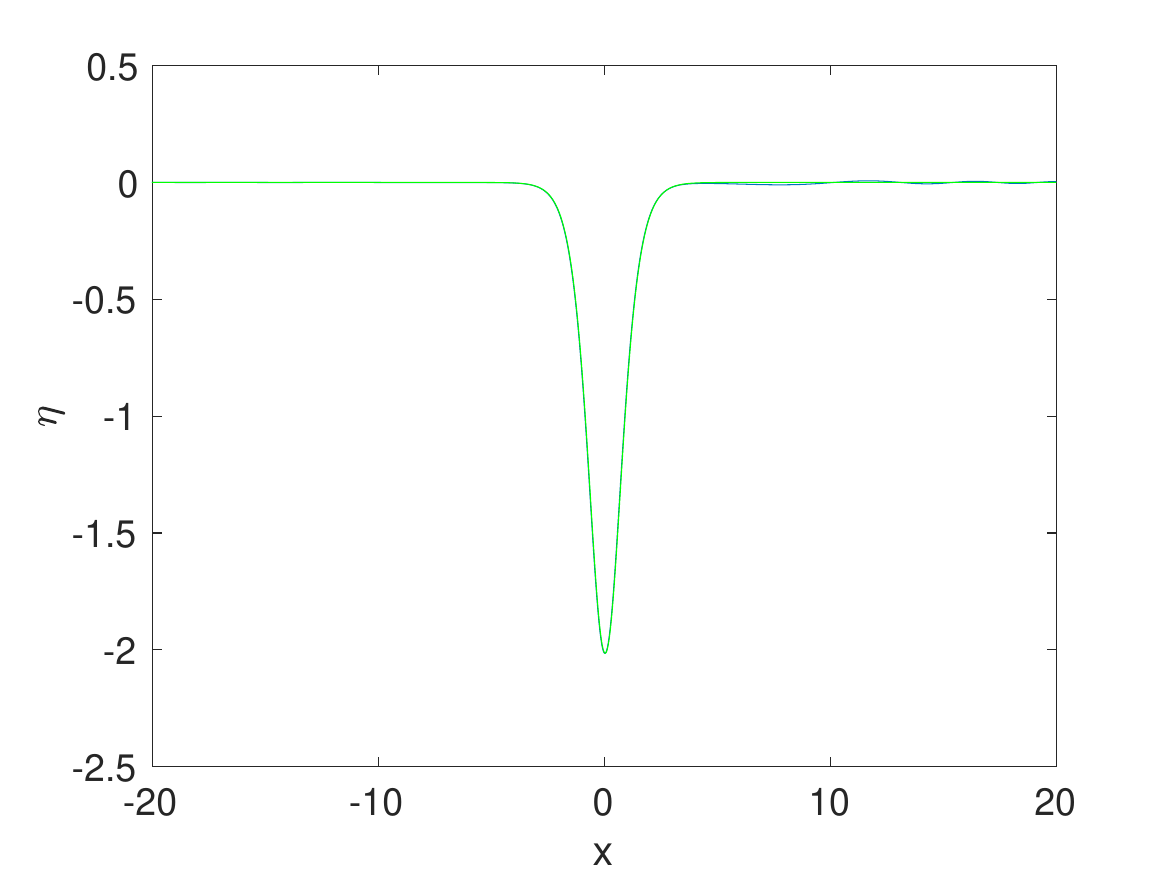}\\
  \includegraphics[width=0.49\textwidth]{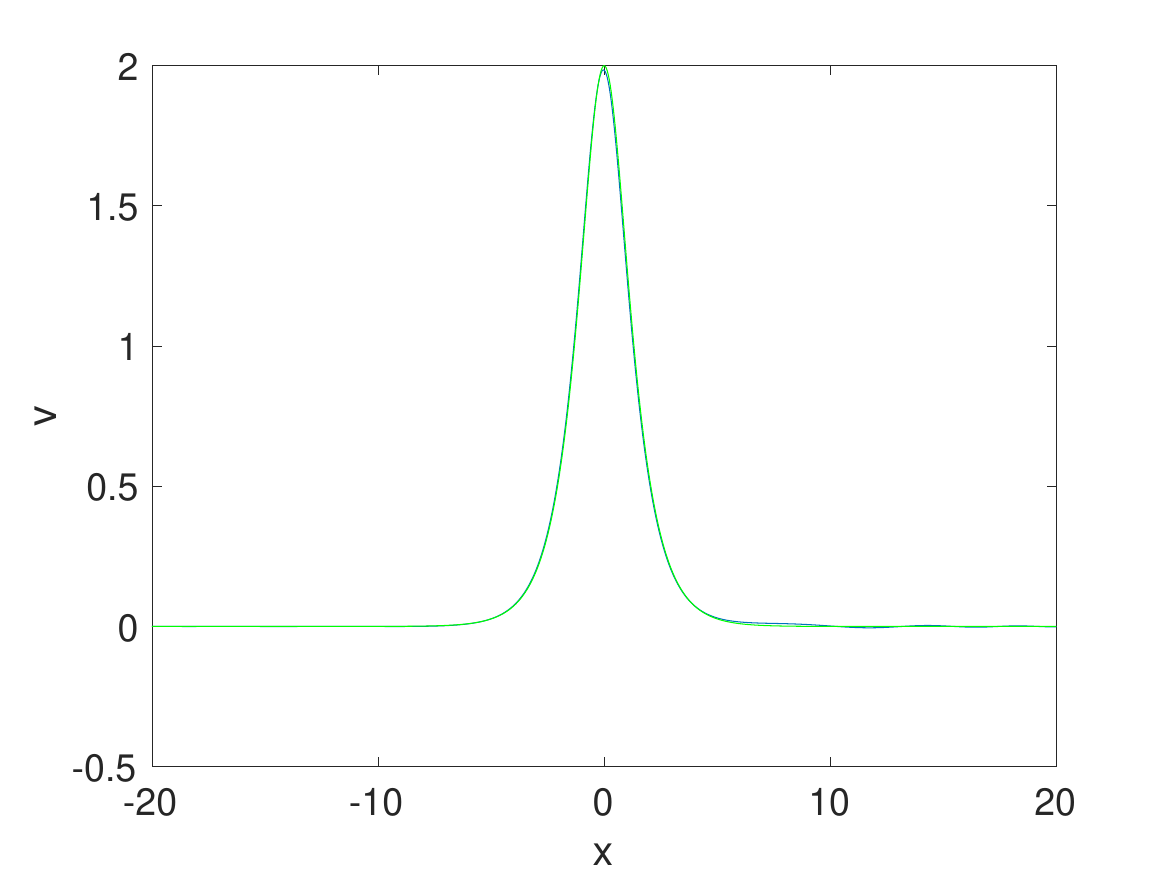}
 \includegraphics[width=0.49\textwidth]{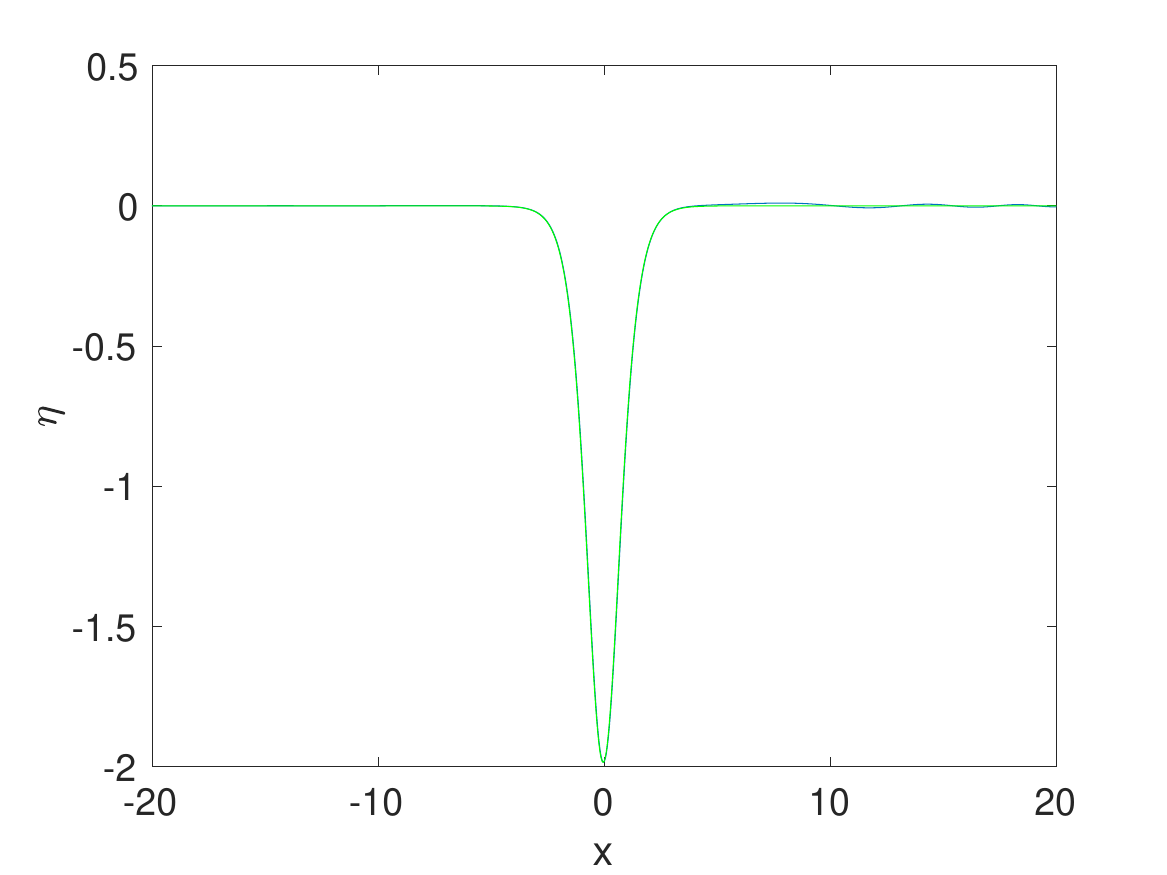}
\caption{The solution to the KBK system for the perturbed  stationary 
solution of the form (\ref{pert}) with $\lambda=1$ and 
 $\mu=1.01$ in the upper row and $\mu=0.99$ in the lower row 
 at the final time $t=5$ in blue and a 
 fitted stationary solution in green, on the left $v$, on the right $\eta$.  }
 \label{KBKstat_lat5}
\end{figure}

For $C=0$, the solutions (\ref{soliton}) become stationary. If we 
apply the same perturbations (\ref{pert}) to this solution, we find 
that the stationary solution is also stable. For $\mu=1$ in 
the initial data (\ref{pert}), we get the solutions shown in 
Fig.~\ref{KBKstat_lat5}.

\begin{figure}[htb!]
 \includegraphics[width=0.49\textwidth]{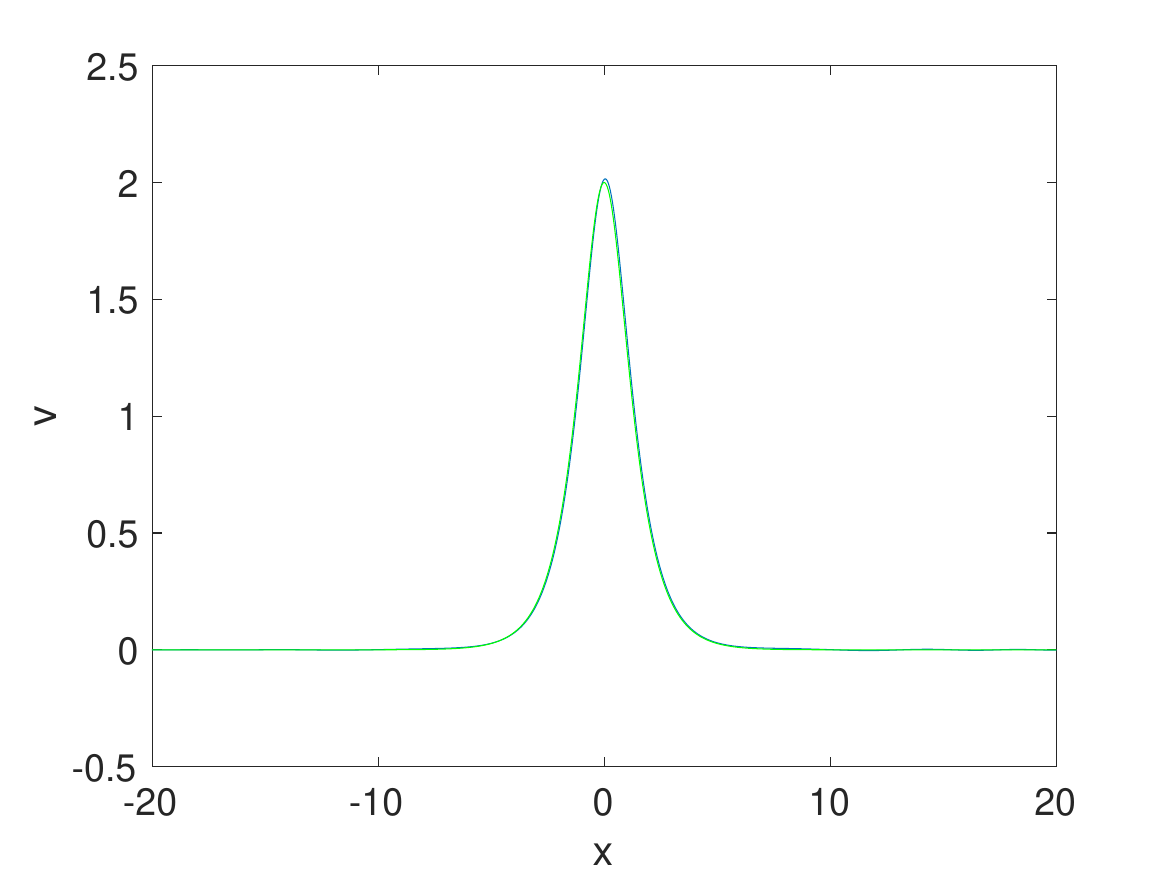}
 \includegraphics[width=0.49\textwidth]{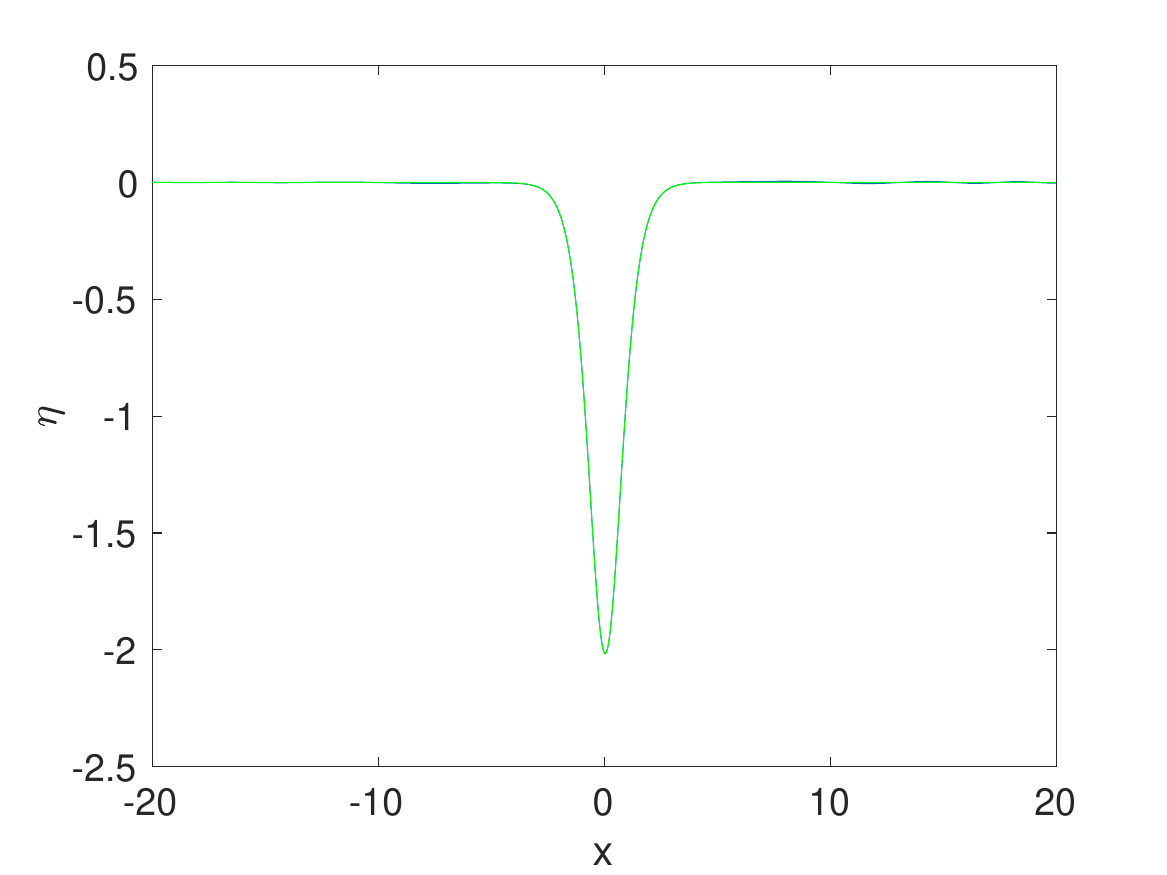}\\
  \includegraphics[width=0.49\textwidth]{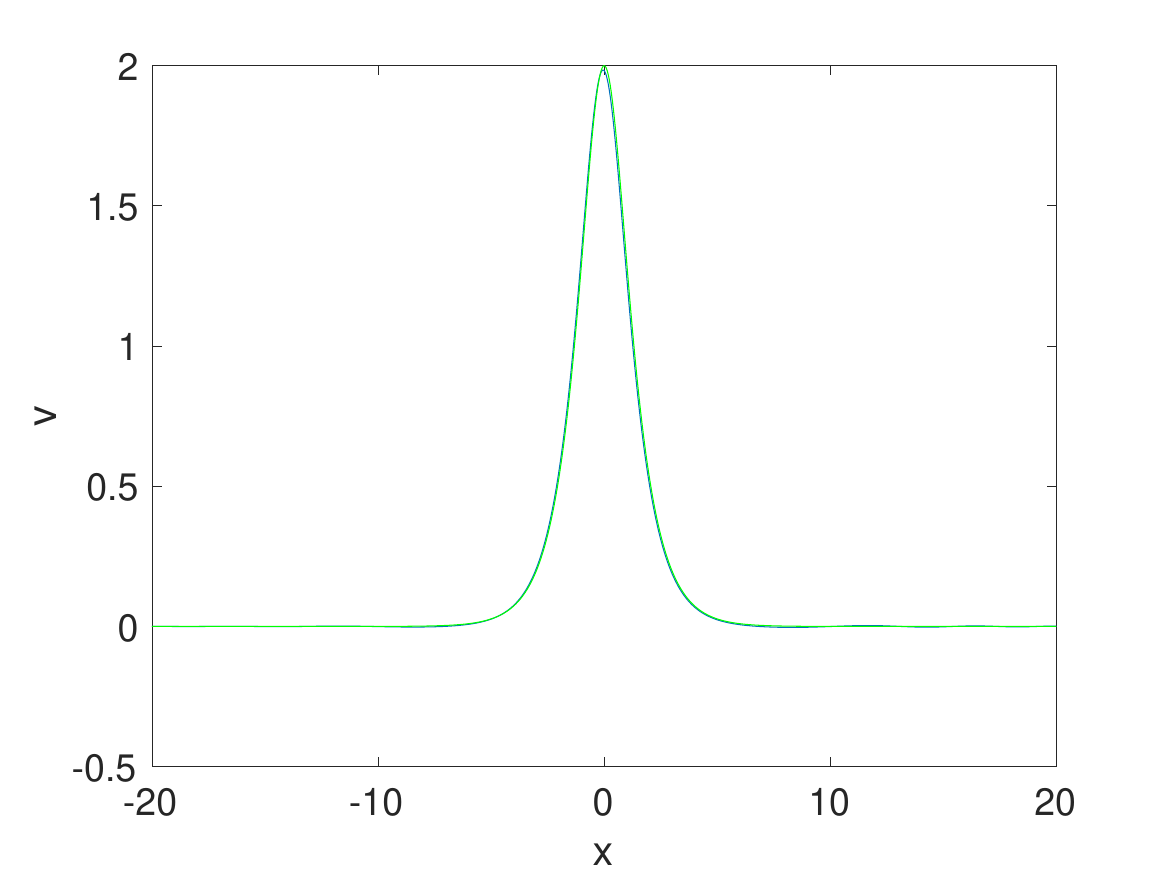}
 \includegraphics[width=0.49\textwidth]{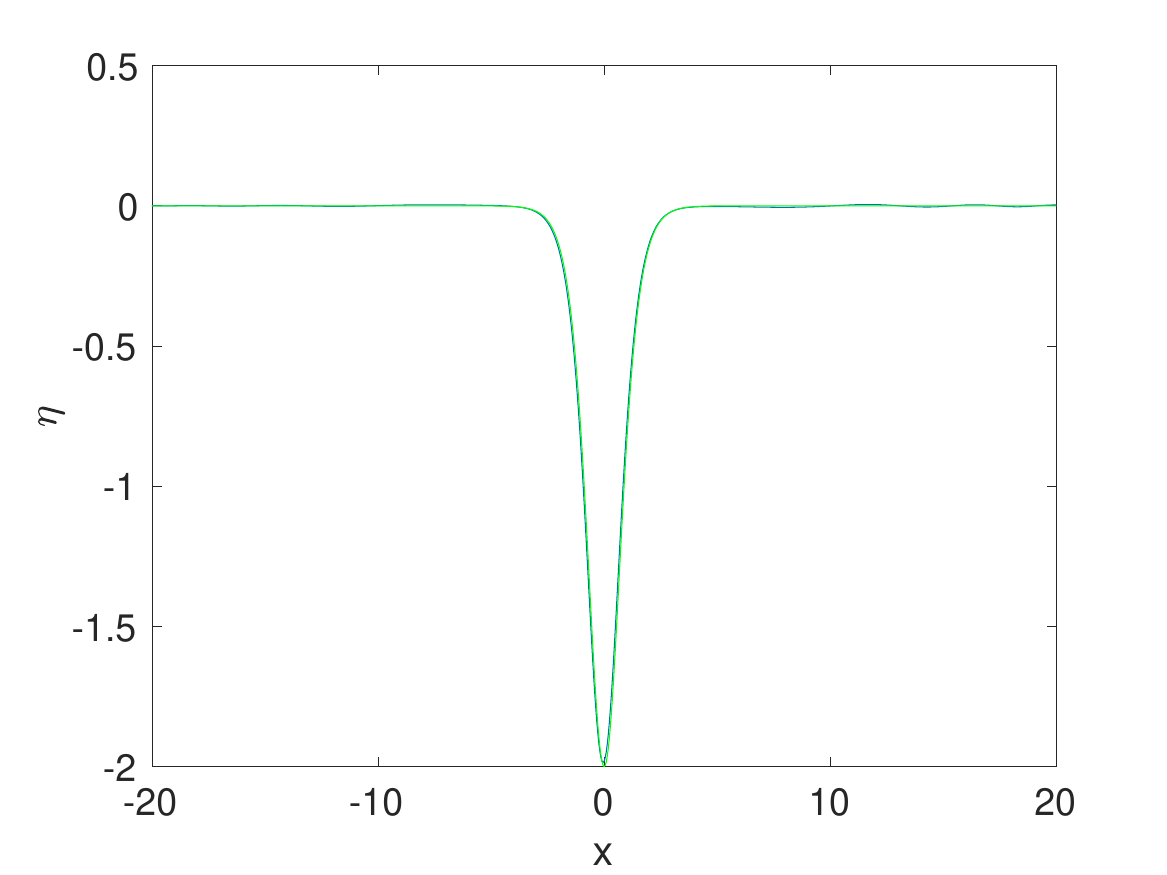}
\caption{The solution to the KBK system for the perturbed  stationary 
solution of the form (\ref{pert}) with $\lambda=1$ and 
$\mu=1.01$ in the upper row and $\mu=0.99$ in the lower row 
 fitted stationary solution in green, on the left $v$, on the right $\eta$.  }
 \label{KBKstat_mut5}
\end{figure}

The situation is similar for $\lambda=1$ as shown in Fig.~\ref{KBKstat_mut5}. 

\section{Localised initial data}
In this section we will study KBK solutions for localised initial 
data. The results confirm the applicability of the soliton resolution 
conjecture to the KBK system, that the long time behavior of 
solutions is given by solitons plus radiation.

\begin{figure}[htb!]
 \includegraphics[width=0.7\textwidth]{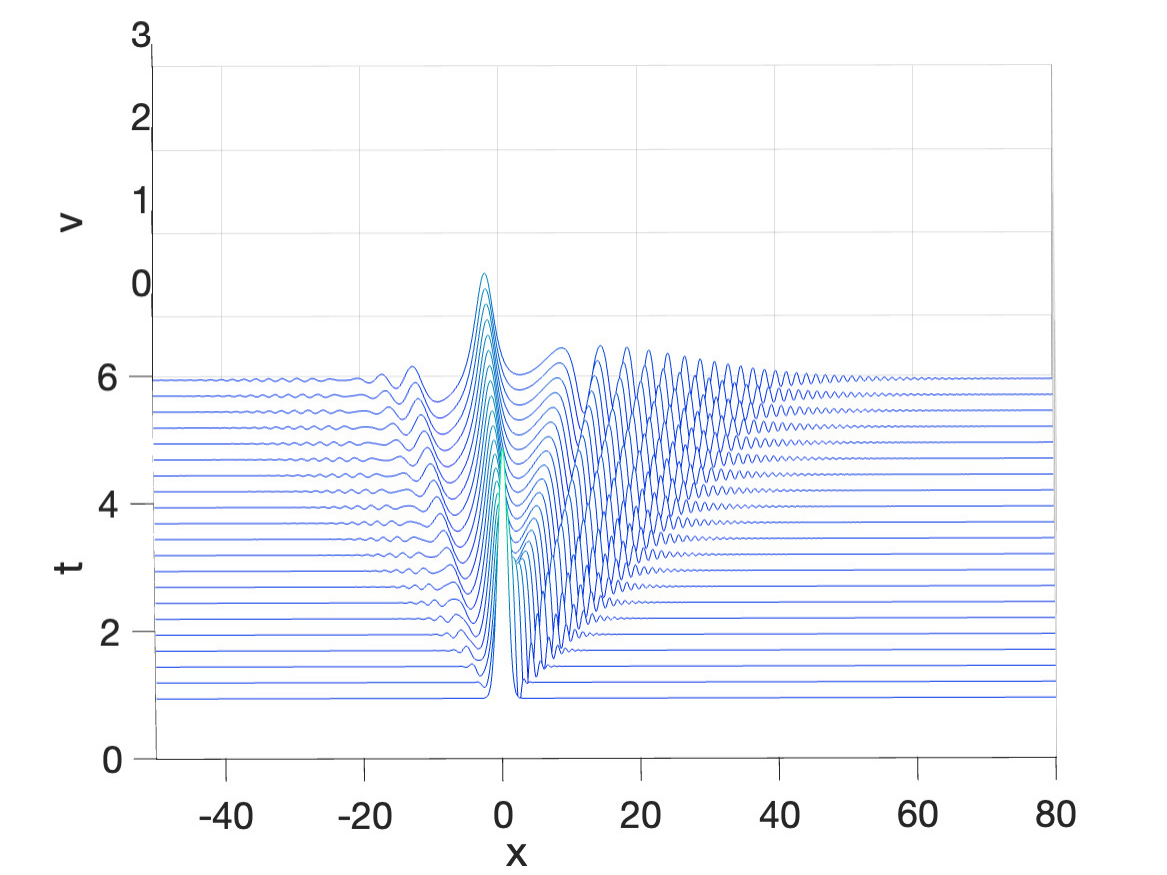}
 \caption{Solution $v$ to the KBK system for the initial data 
 (\ref{vgauss}). }
 \label{v3gausswaterv}
\end{figure}

In this section we always use $N=2^{12}$ DFT modes for $x\in 
30[-\pi,\pi]$ and $N_{t}=4000$ 
time steps. The DFT coefficients decrease to machine precision in all 
examples, and the relative conservation of the energy is of the order 
of $10^{-10}$ and better during the computations. First we consider 
initial data of the form 
\begin{equation}
	\eta(x,0)=0,\quad v(x,0) = 3\exp(-x^{2})
	\label{vgauss}.
\end{equation}
The resulting solution for $v$ can be seen in 
Fig.~\ref{v3gausswaterv}. There are strong oscillations propagating 
to the right and smaller ones propagating to the left. And there is a 
solitary wave traveling towards $-\infty$ to be discussed in more 
detail below. 

\begin{figure}[htb!]
 \includegraphics[width=0.7\textwidth]{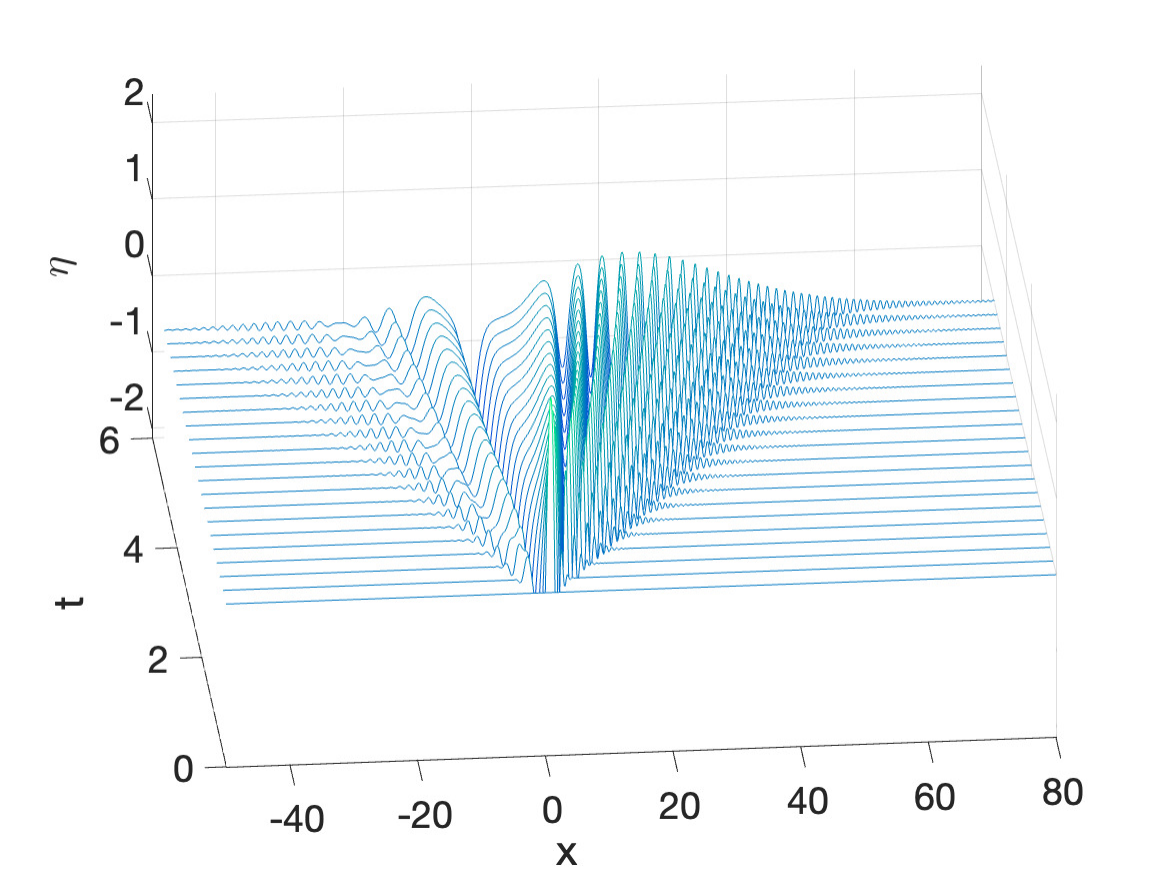}
 \caption{Solution $\eta$ to the KBK system for the initial data 
 (\ref{vgauss}). }
 \label{v3gausswatereta}
\end{figure}

The corresponding solution $\eta$ can be seen in 
Fig.~\ref{v3gausswatereta}. There is similar radiation as in the 
solution $v$, and again a solitary structure traveling to the left. 

\begin{figure}[htb!]
 \includegraphics[width=0.49\textwidth]{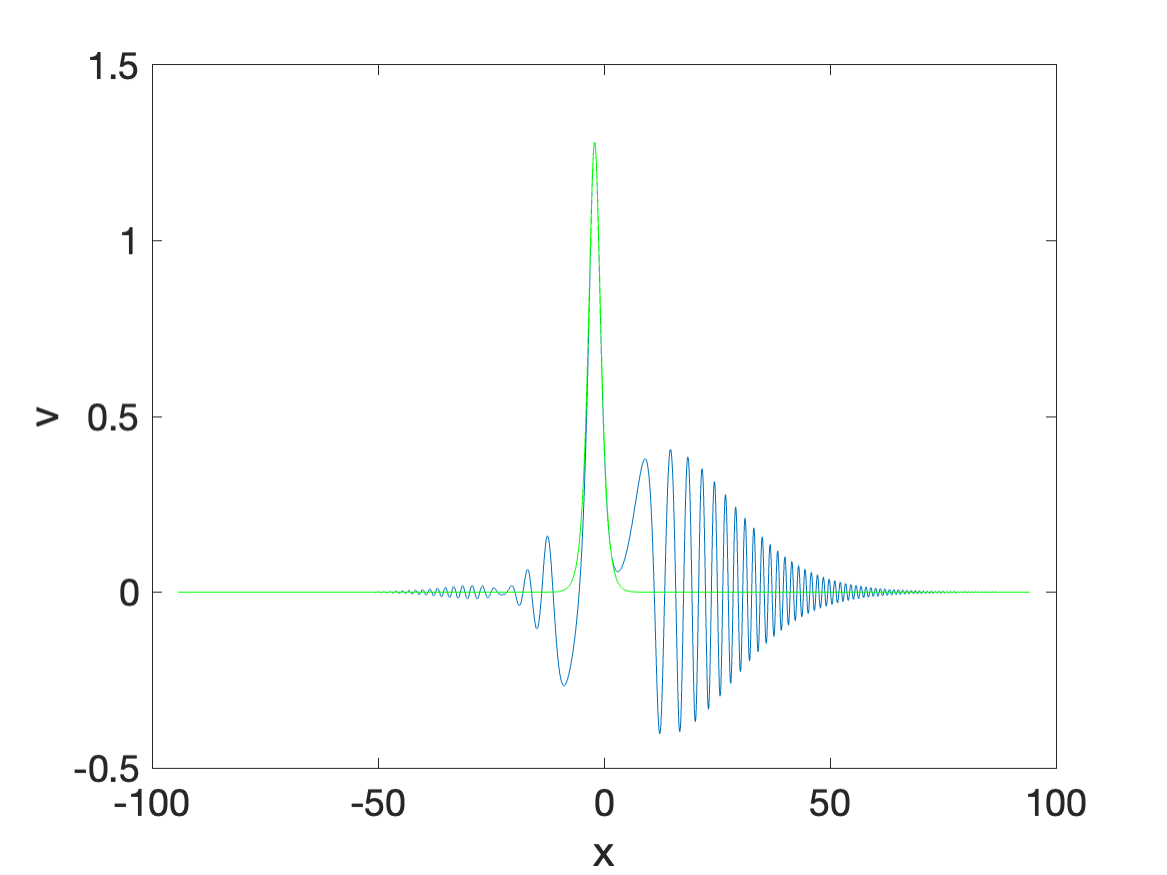}
 \includegraphics[width=0.49\textwidth]{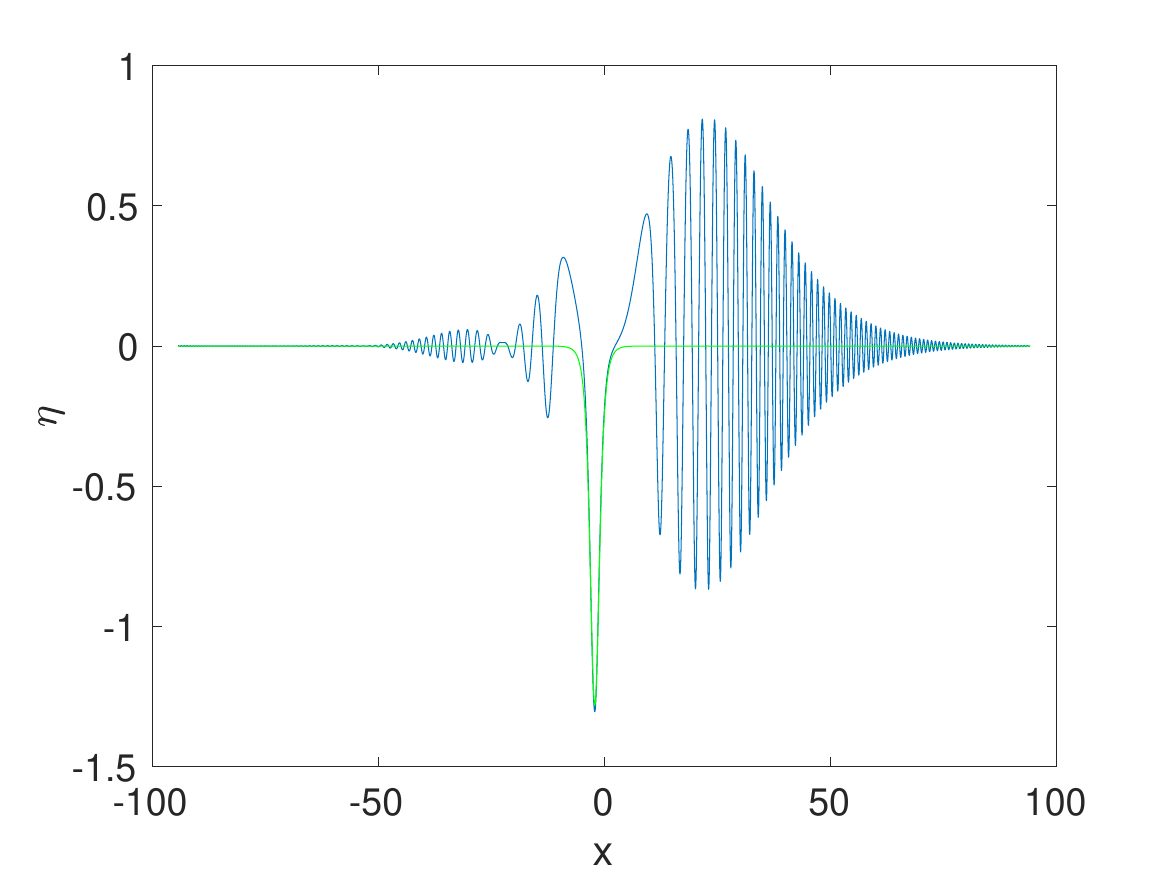}
 \caption{The solutions shown in Fig.~\ref{v3gausswaterv} and 
 Fig.~\ref{v3gausswatereta} at the final time $t=5$ in blue and a 
 fitted soliton in green.  }
 \label{KBKv3gausst5solfit}
\end{figure}

To check whether the solitary structure near the origin is in fact 
evolving into a soliton, we fit it to the soliton (\ref{soliton}) as 
described in the previous section. 
We show the result of this fitting for $v$ on left of 
Fig.~\ref{KBKv3gausst5solfit}. The numerical solution is given in 
blue, the fitted soliton in green. The corresponding plot for $\eta$ is 
shown on the right of the same figure. It can be seen that the 
agreement is already very good though there is still a considerable 
amount of radiation in the vicinity of the soliton. 

\begin{figure}[htb!]
 \includegraphics[width=0.7\textwidth]{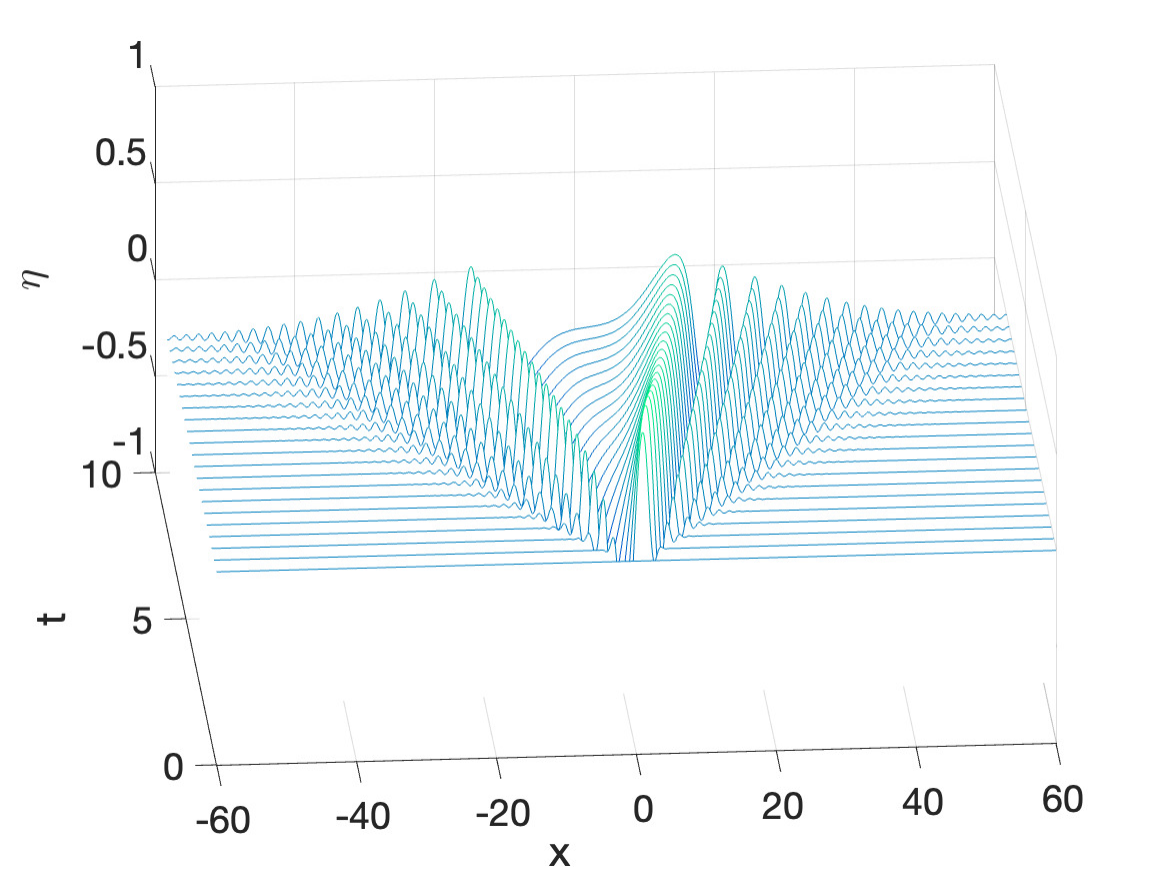}
 \caption{Solution $v$ to the KBK system for the initial data 
 (\ref{etagauss}) with $A=3$. }
 \label{KBKeta3gausswaterv}
\end{figure}

Note that there is no criterion known which initial data lead to 
which number of KBK solitons for large times. If we consider for 
instance the 
initial data 
\begin{equation}
	\eta(x,0)=A\exp(-x^{2}),\quad v(x,0) = 0,\quad A\in\mathbb{R}
	\label{etagauss},
\end{equation}
we get with $A=3$  for $v$ the solution shown in Fig.~\ref{KBKeta3gausswaterv}. In 
this case there is no indication of solitons. 

\begin{figure}[htb!]
 \includegraphics[width=0.7\textwidth]{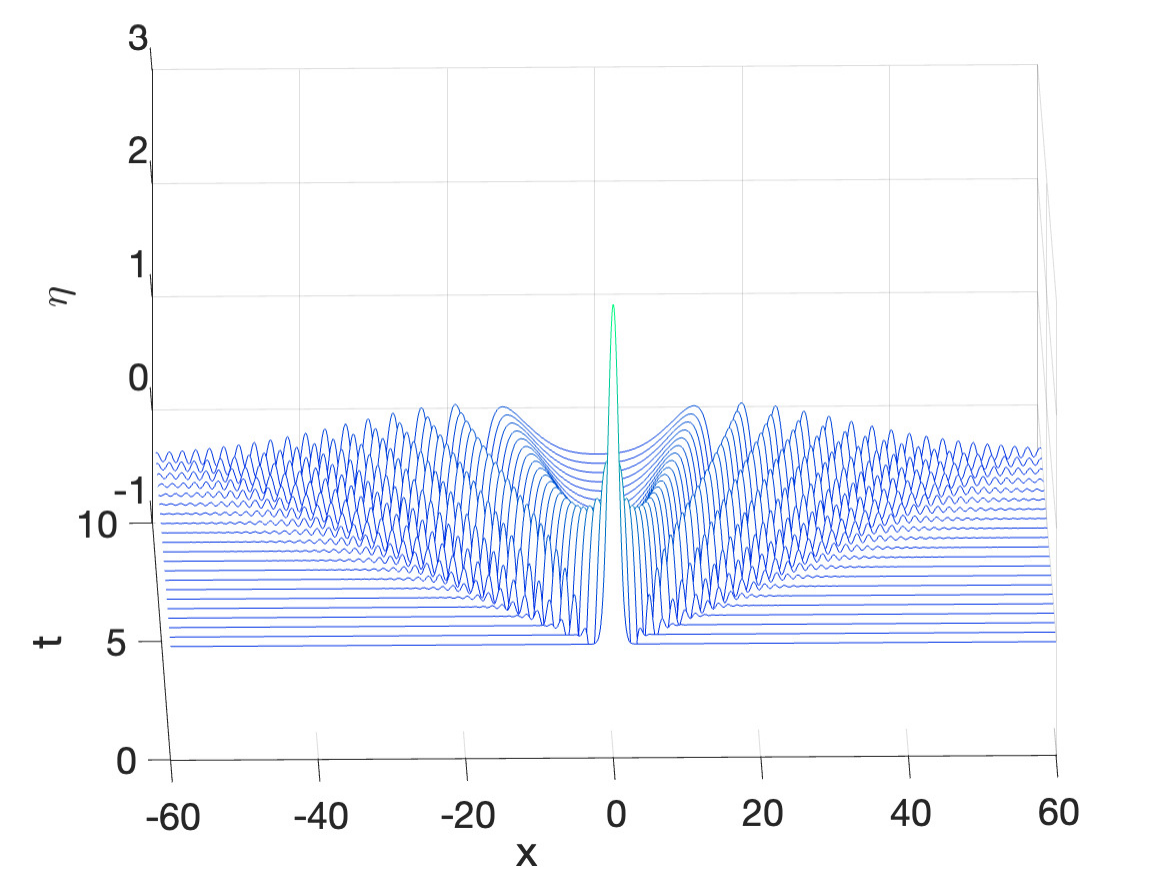}
 \caption{Solution $\eta$ to the KBK system for the initial data 
 (\ref{etagauss}) with $A=3$. }
 \label{KBKeta3gausswatereta}
\end{figure}

The corresponding solution $\eta$ is shown in 
Fig.~\ref{KBKeta3gausswatereta}. It appears that the initial data are 
simply dispersed in this case. Note that the solution looks 
qualitatively similar for larger values for $A$, for instance $A=20$, there is no 
indication of solitons in this class of initial data. 

\begin{figure}[htb!]
 \includegraphics[width=0.7\textwidth]{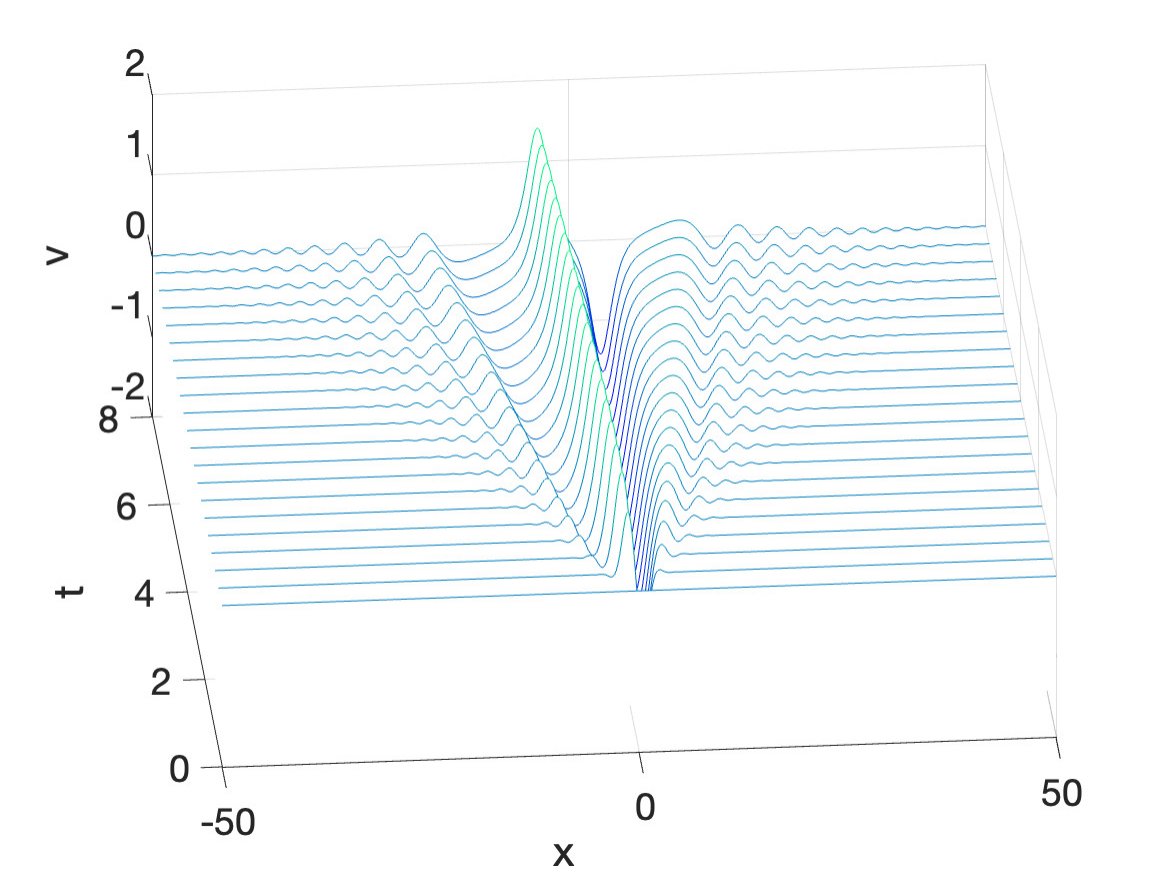}
 \caption{Solution $v$ to the KBK system for the initial data 
 (\ref{etagauss}) with $A=-3$. }
 \label{KBKetam3gausswaterv}
\end{figure}

However, the situation is different for negative $A$ in 
(\ref{etagauss}) as can be seen for $A=-3$ in 
Fig.~\ref{KBKetam3gausswaterv}. There appears to be a solitary 
structure near the origin.

The corresponding solution $\eta$ can be seen in Fig.~\ref{KBKeta3gausswatereta}.
 \begin{figure}[htb!]
 \includegraphics[width=0.7\textwidth]{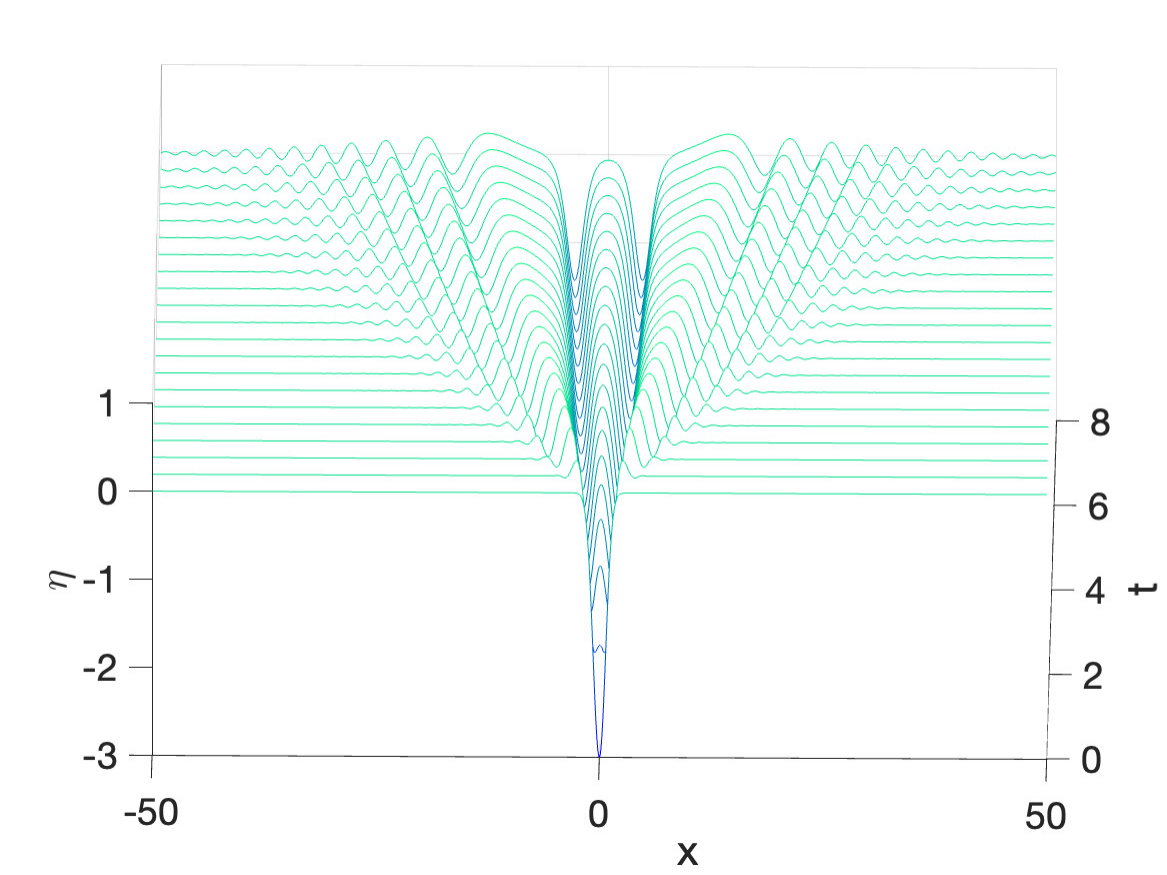}
 \caption{Solution $\eta$ to the KBK system for the initial data 
 (\ref{etagauss}) with $A=-3$. }
 \label{KBKetam3gausswatereta}
\end{figure}

\begin{figure}[htb!]
 \includegraphics[width=0.49\textwidth]{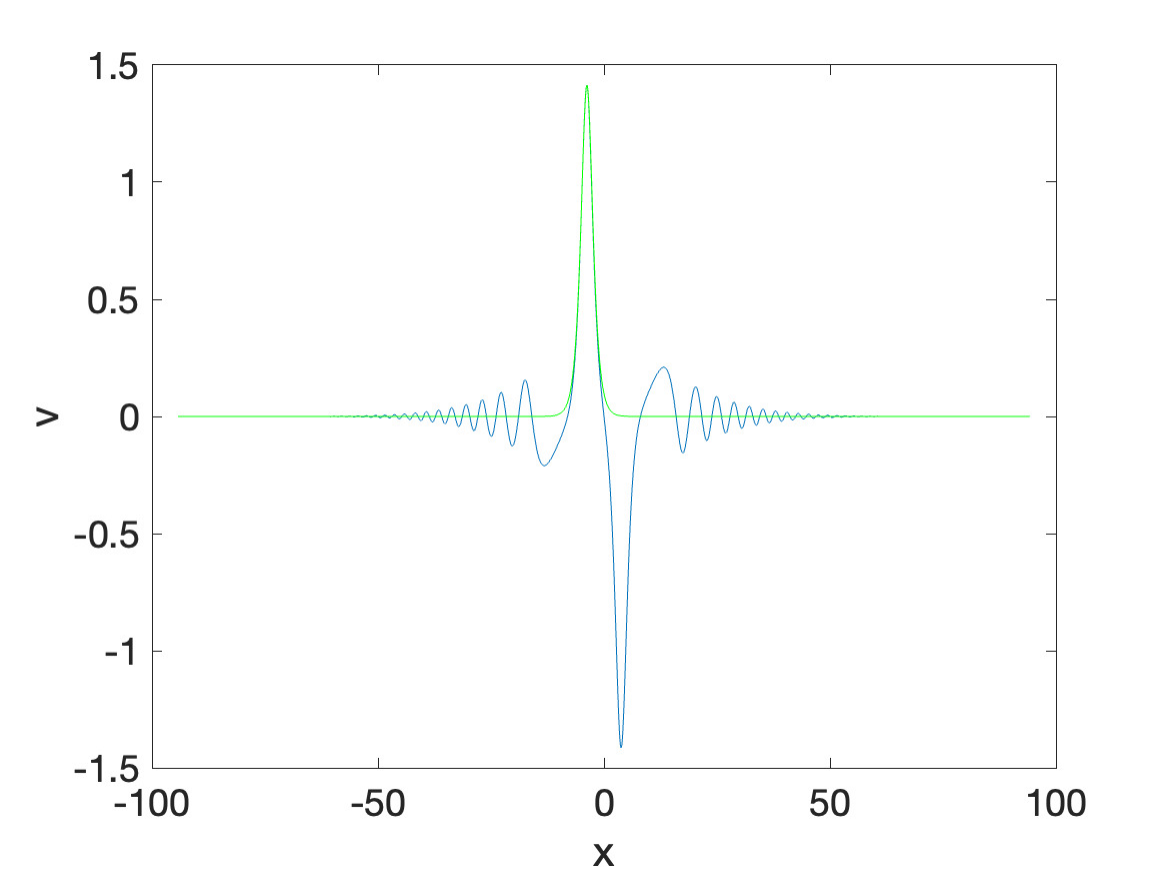}
 \includegraphics[width=0.49\textwidth]{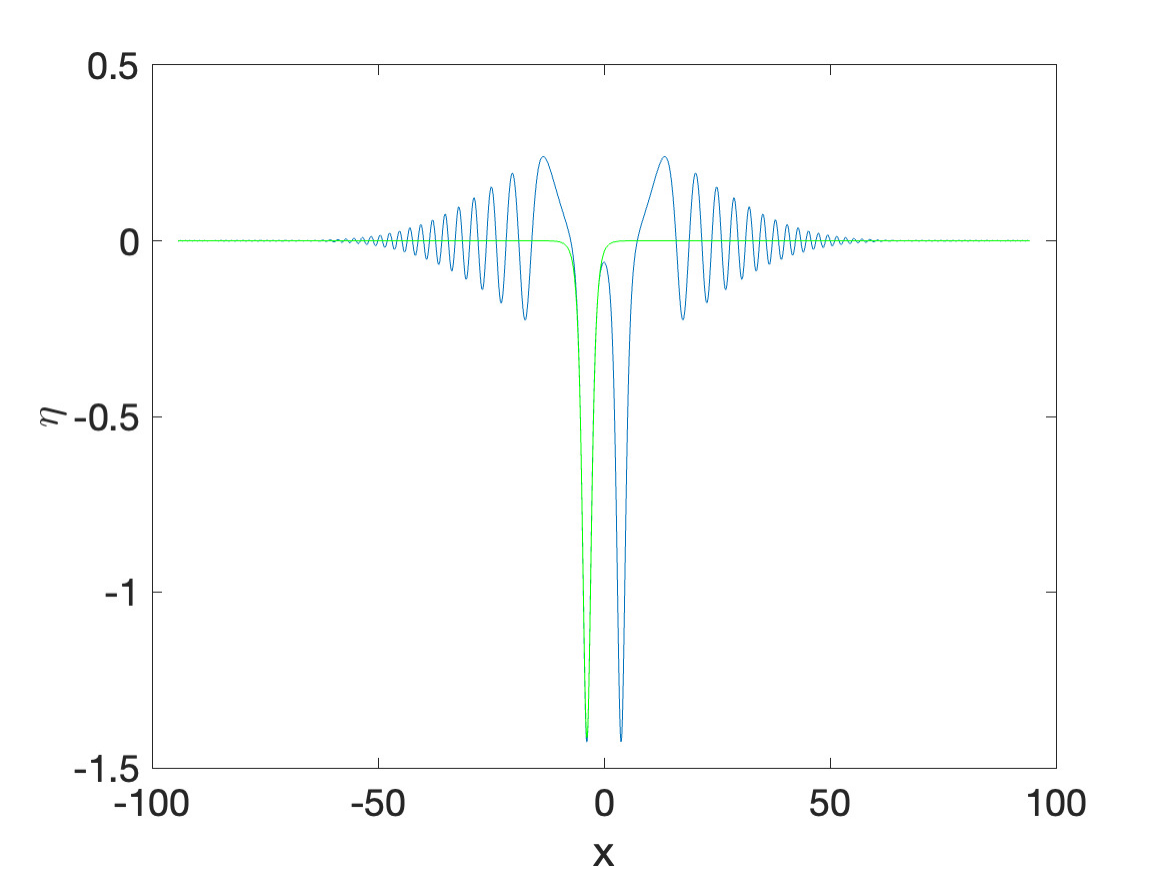}
 \caption{The solutions shown in Fig.~\ref{KBKetam3gausswaterv} and 
 Fig.~\ref{KBKetam3gausswatereta} at the final time $t=8$ in blue and a 
 fitted soliton in green.  }
 \label{KBKetam3gausst8solfit}
\end{figure}

As in Fig.~\ref{KBKv3gausst5solfit}, we fit the left maximum to the 
soliton (\ref{soliton}). We show the solution for $v$ at $t=8$ with 
the fitted soliton in green on the left of 
Fig.~\ref{KBKetam3gausst8solfit}. The corresponding plot for $\eta$ 
can be seen on the right of the same figure. There is clearly a 
second soliton moving to  the right in this case for symmetry 
reasons. 

Note that these initial data do not 
satisfy the non-cavitation condition, but that there is no indication 
of a blow-up in this case. The solution appears to be smooth for all 
times.

\section{Dispersive shock waves}
Dispersive shock waves (DSWs) are zones of rapid modulated oscillations in 
solutions to dispersive nonlinear PDEs in the vicinity of shocks of 
the corresponding dispersionless system, here the Saint-Venant 
system. A possible way to study such zones for the KBK system is to 
consider solutions on $x$ scales of order $1/\varepsilon$ for 
$\varepsilon\ll1$ on time scales of order $1/\varepsilon$. This can be 
conveniently done by rescaling $x$ and $t$ by a factor $1/\varepsilon$. 
This leads for (\ref{KBK}) (we use the same notation as before)
\begin{equation}
    \label{KBKe}
    \left\lbrace
    \begin{array}{l}
    \eta_t+v_x+(\eta v)_x+\varepsilon^{2} v_{xxx}=0 \\
    v_t+\eta_x+vv_x=0.
\end{array}\right.
    \end{equation}
    
In the formal limit $\varepsilon\to0$, this leads to the Saint-Venant 
system
\begin{equation}
    \label{SV}
    \left\lbrace
    \begin{array}{l}
    \eta_t+v_x+(\eta v)_x=0 \\
    v_t+\eta_x+vv_x=0.
\end{array}\right.
    \end{equation}
the solutions of which will have shocks in finite time for hump-like 
initial data. The system (\ref{KBKe}) can be seen as a dispersive 
regularization of the system (\ref{SV}). For the same initial data 
leading to shocks for the Saint-Venant system, DSWs are expected in 
the vicinity of the former for the system (\ref{KBKe}). 

\begin{figure}[htb!]
 \includegraphics[width=0.7\textwidth]{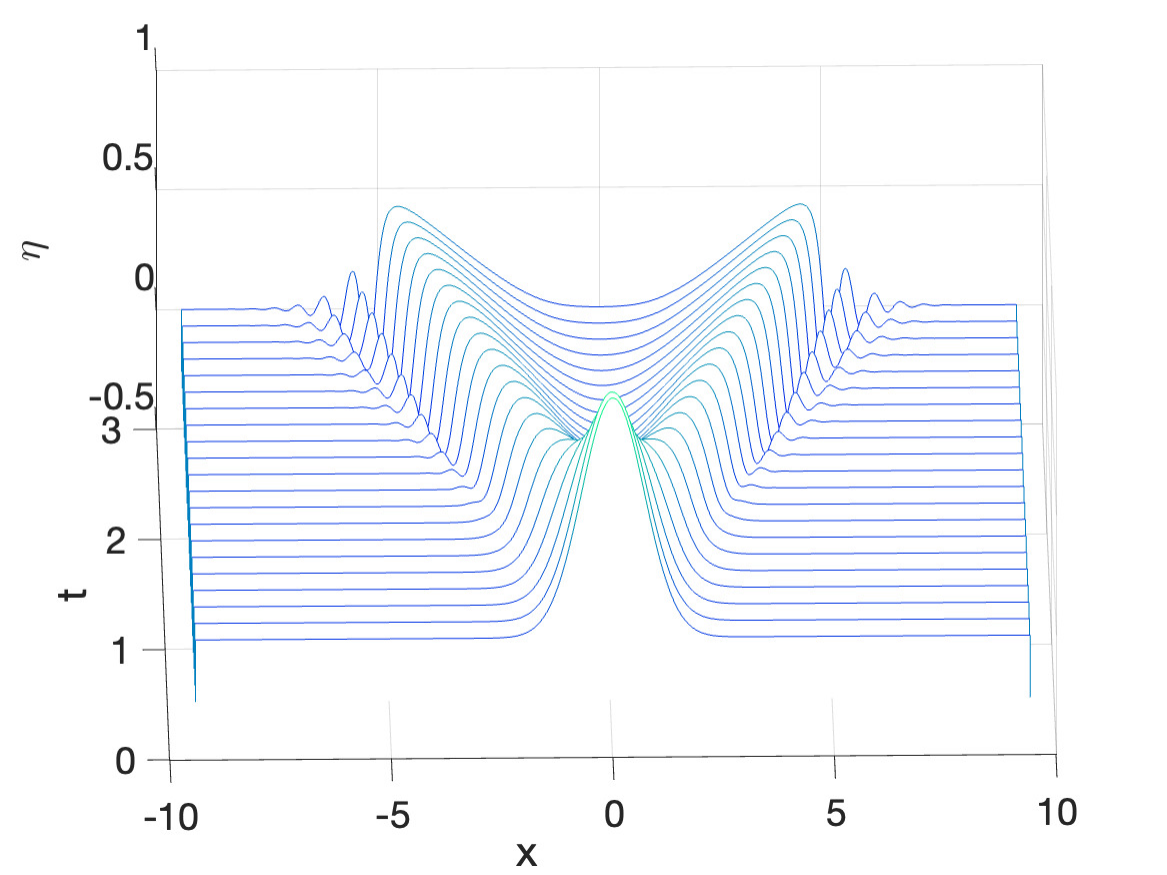}
 \caption{Solution $\eta$ to the KBK system (\ref{KBKe}) for 
 $\varepsilon=0.1$ and for the initial data 
 $\eta(x,0)=\exp(-x^{2})$, 
$v(x,0)=0$. }
 \label{KBKetagauss1e2watereta}
\end{figure}

We use $N=2^{14}$ DFT modes for $x\in 3[-\pi,\pi]$ and $10^{4}$ time 
steps for $t\leq 3$ for the initial data $\eta(x,0)=\exp(-x^{2})$, 
$v(x,0)=0$. The solution $\eta$ for $\varepsilon=0.1$ can be seen in 
Fig.~\ref{KBKetagauss1e2watereta}. The initial hump splits into two 
humps developing strong gradients at the outer edges where 
oscillations can be observed. 

\begin{figure}[htb!]
 \includegraphics[width=0.7\textwidth]{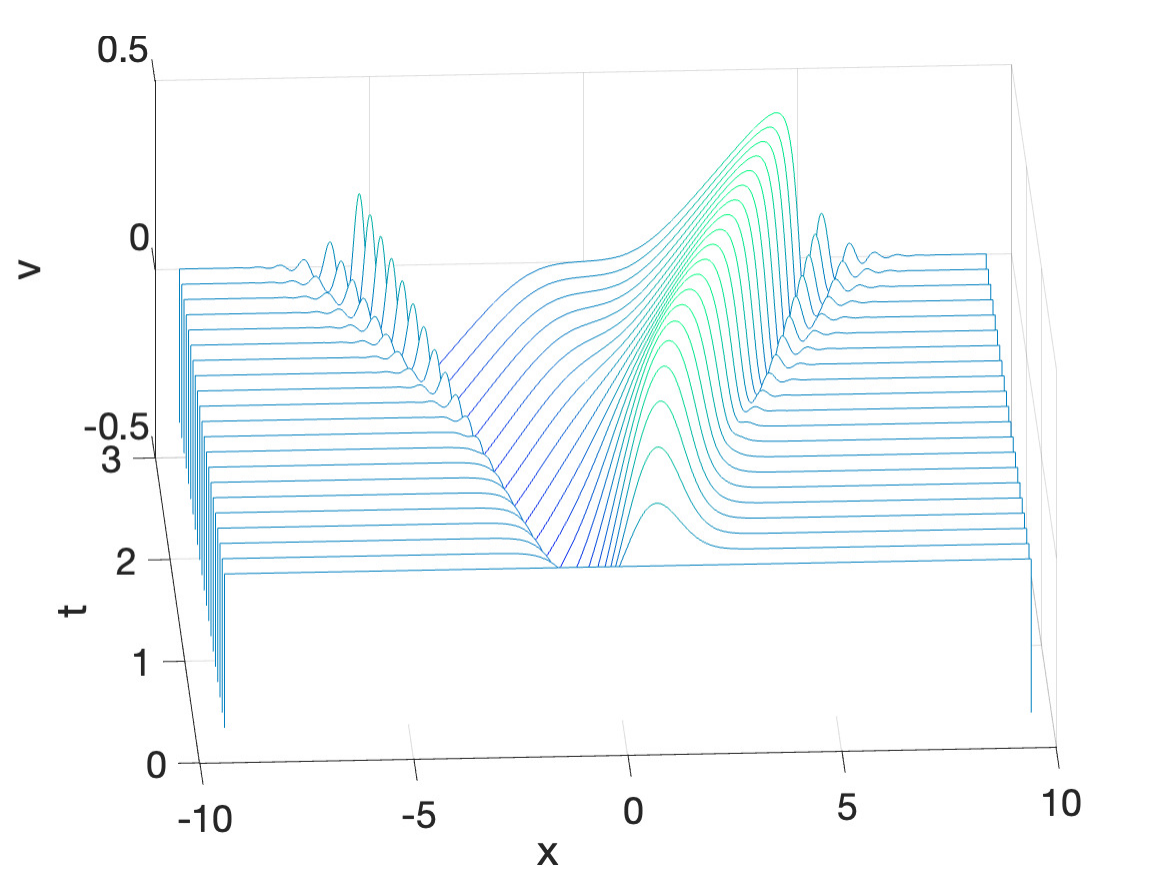}
 \caption{Solution $\eta$ to the KBK system (\ref{KBKe}) for 
 $\varepsilon=0.1$ and for the initial data 
 $\eta(x,0)=\exp(-x^{2})$, 
$v(x,0)=0$. }
 \label{KBKetagauss1e2waterv}
\end{figure}

The  corresponding solution $v$ is shown in 
Fig.~\ref{KBKetagauss1e2waterv}. The solution is odd in $x$, but 
shows oscillations at the same $x$-values as $\eta$. 

\begin{figure}[htb!]
 \includegraphics[width=0.49\textwidth]{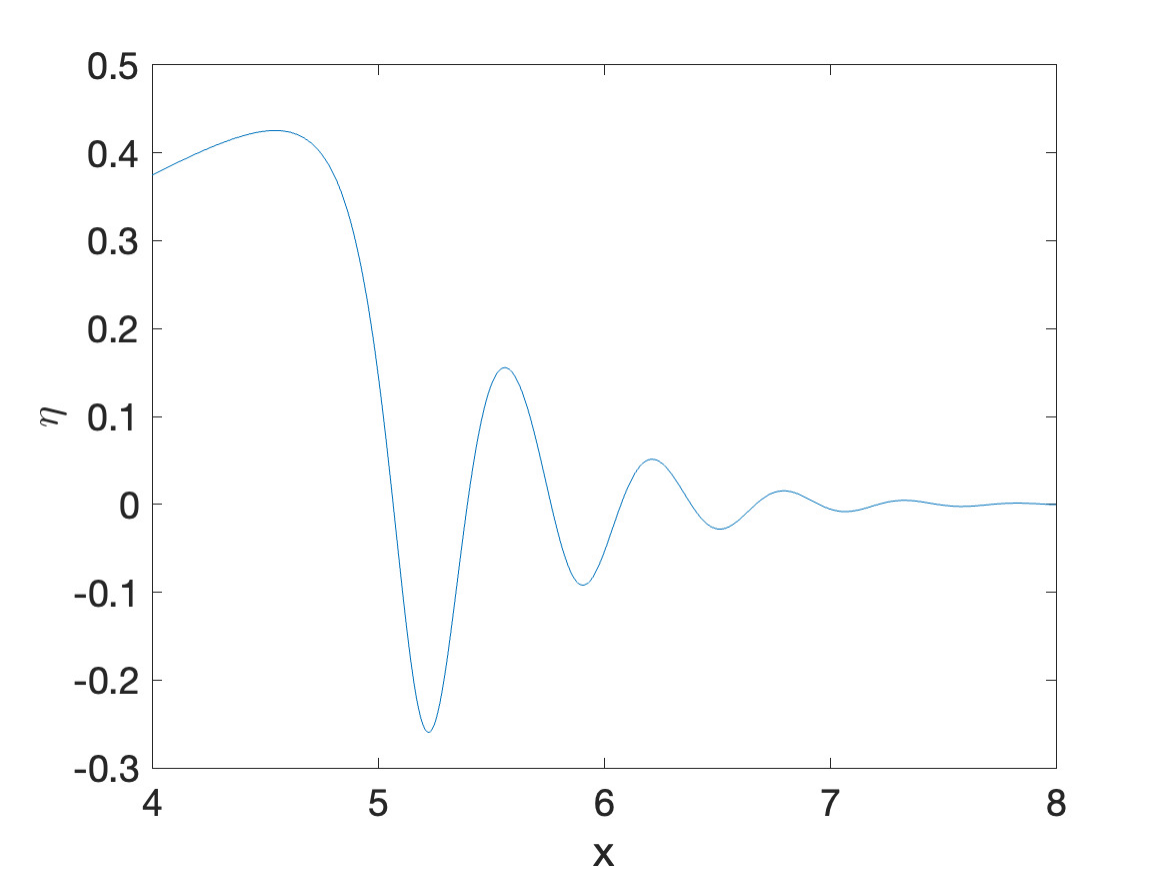}
 \includegraphics[width=0.49\textwidth]{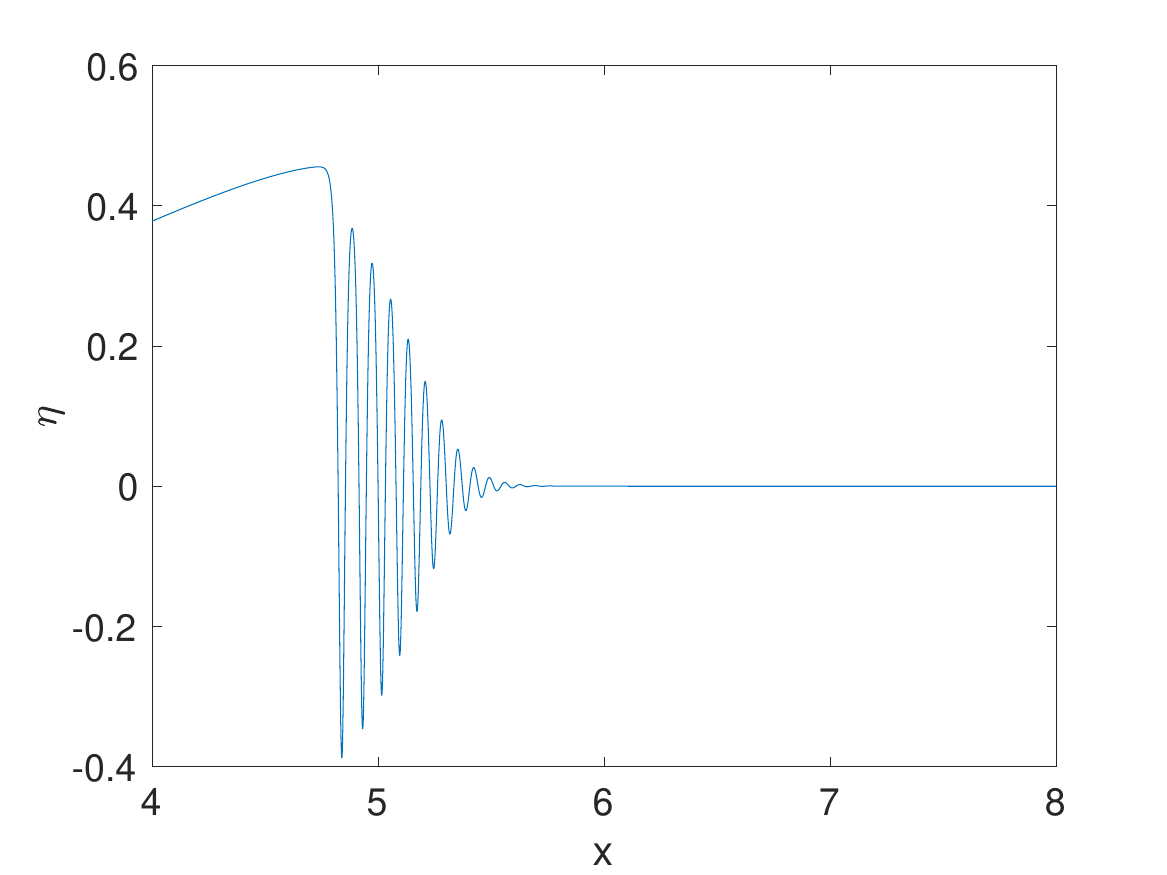}\\
\includegraphics[width=0.49\textwidth]{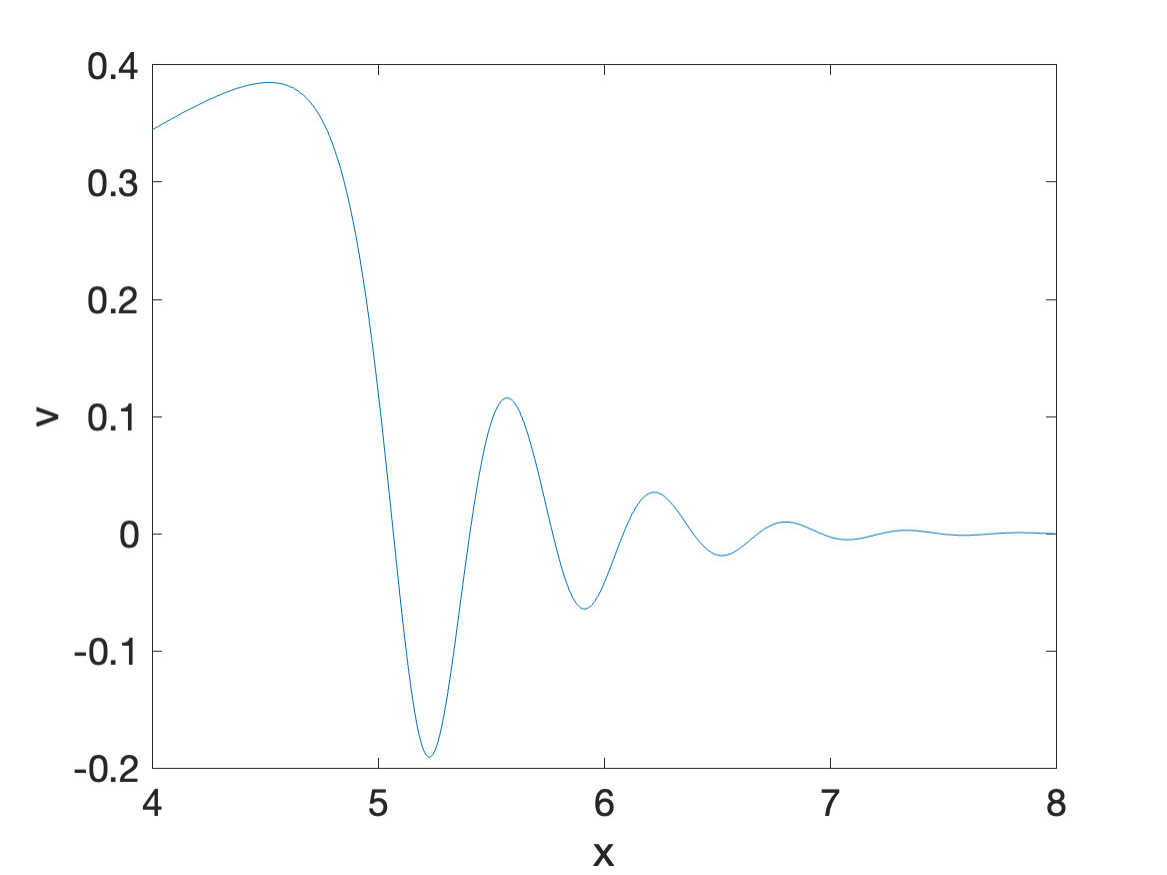}
 \includegraphics[width=0.49\textwidth]{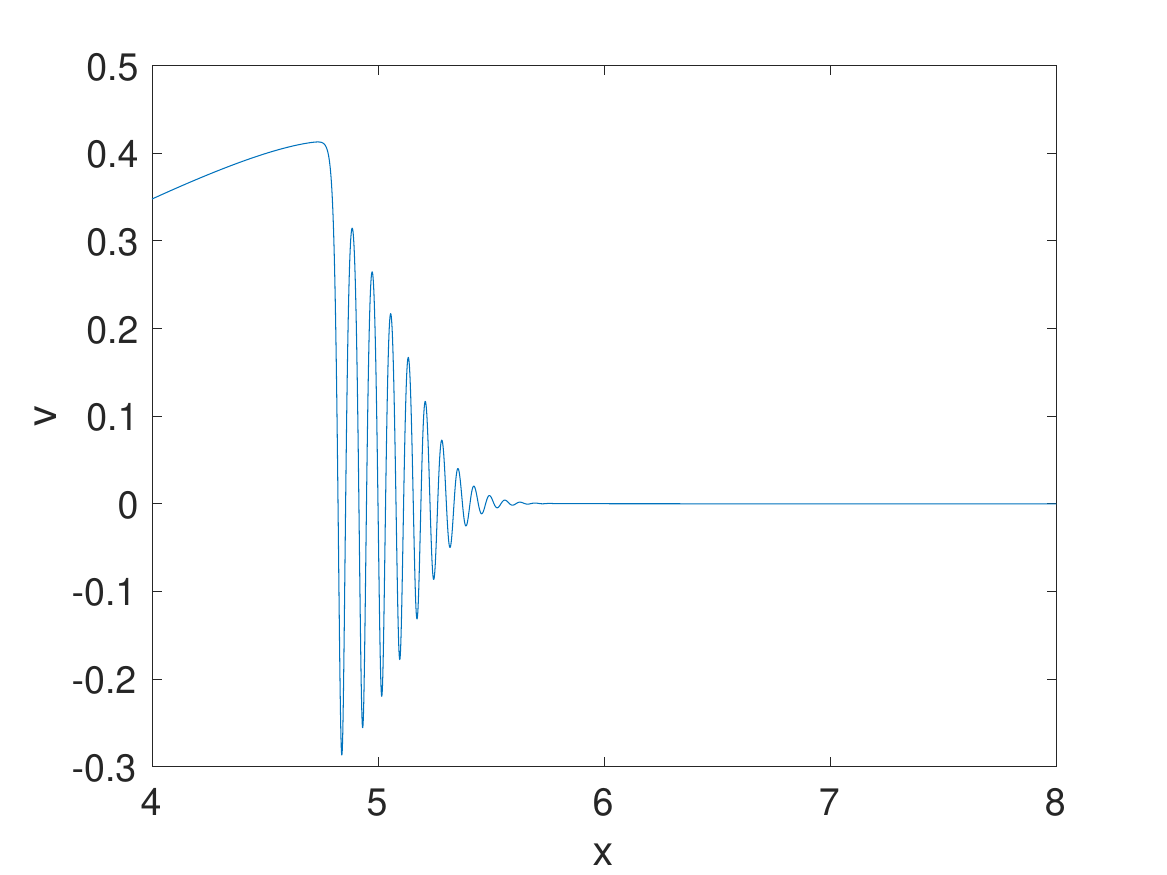}
 \caption{On the left a close-up of the solutions in 
 Fig.~\ref{KBKetagauss1e2watereta} and 
 Fig.~\ref{KBKetagauss1e2waterv} at the final time $t=3$, in the upper row $\eta$, in the 
 lower row $v$; on the right the  oscillatory zone for the same 
 initial data for $\varepsilon=0.01$. }
 \label{KBKetagausszoom}
\end{figure}

We show a close up of the oscillatory zone in 
Fig.~\ref{KBKetagausszoom}. The smaller $\varepsilon$, the more rapid 
the oscillations and the more they are localised to a zone sometimes 
called the Whitham zone. 

\begin{figure}[htb!]
 \includegraphics[width=0.7\textwidth]{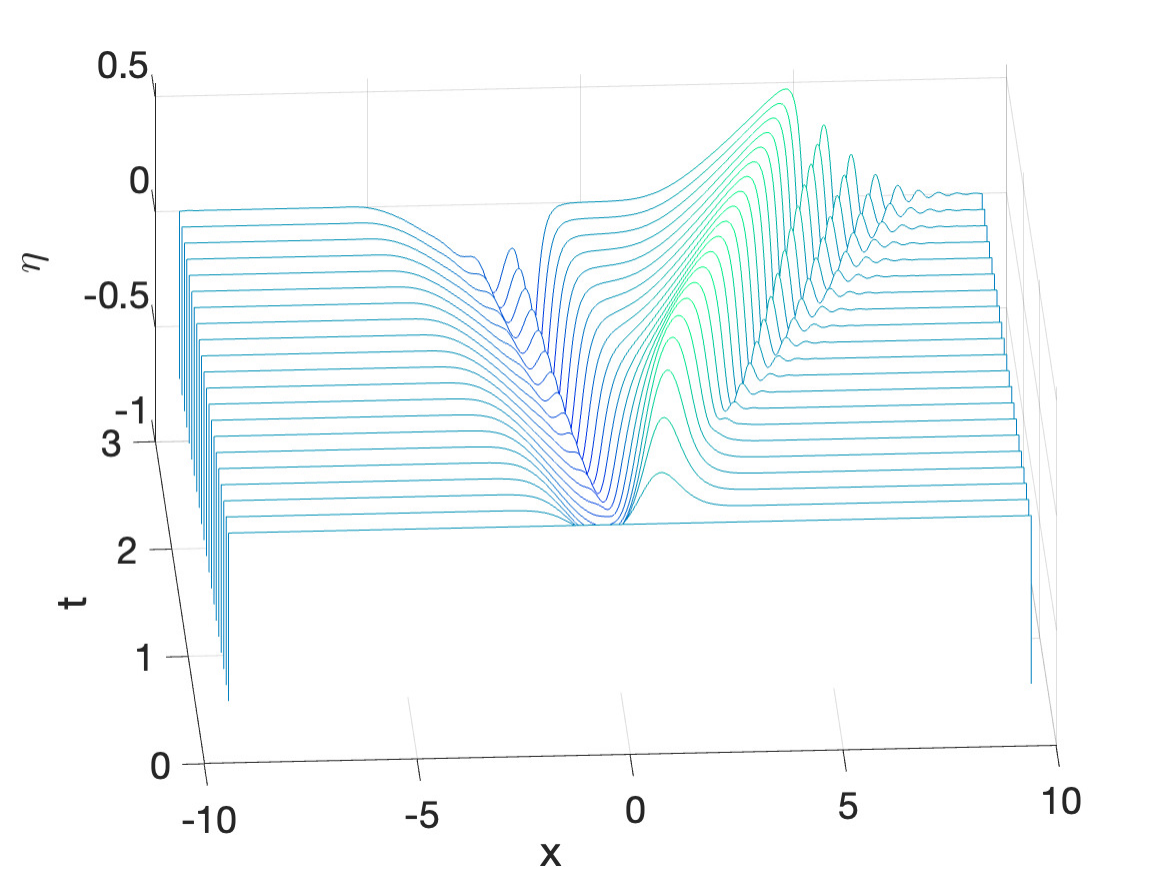}
 \caption{Solution $\eta$ to the KBK system (\ref{KBKe}) for 
 $\varepsilon=0.1$ and for the initial data 
 $v(x,0)=\exp(-x^{2})$, 
$\eta(x,0)=0$. }
 \label{KBKvgauss1e2watereta}
\end{figure}

The situation is similar for other hump-like initial data as 
$\eta(x,0)=0$ and $v(x,0)=\exp(-x^{2})$. The solution $\eta$ for 
$\varepsilon=0.1$ can be seen in Fig.~\ref{KBKvgauss1e2watereta}. The 
solution is not symmetric, but there are oscillatory zones as before 
on both sides of the humps. 

\begin{figure}[htb!]
 \includegraphics[width=0.7\textwidth]{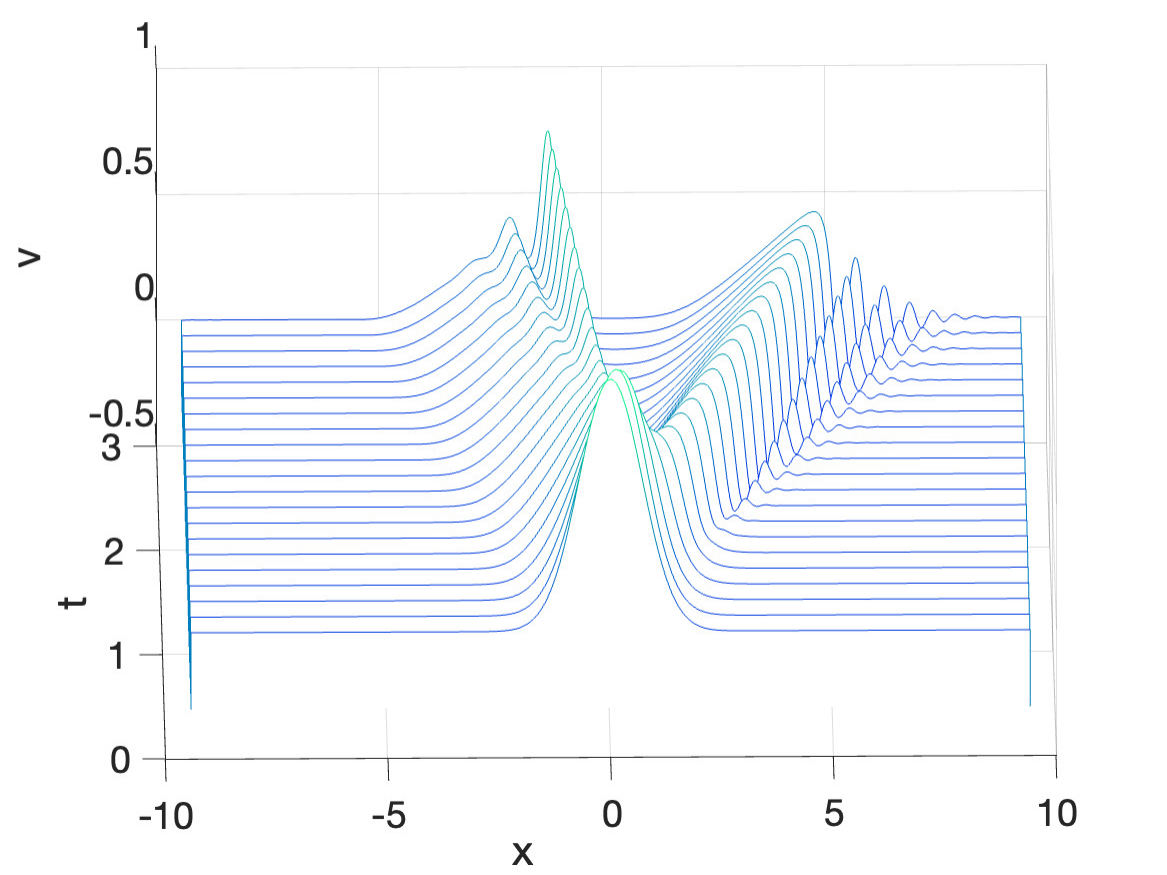}
 \caption{Solution $v$ to the KBK system (\ref{KBKe}) for 
 $\varepsilon=0.1$ and for the initial data 
 $v(x,0)=\exp(-x^{2})$, 
$\eta(x,0)=0$. }
 \label{KBKvgauss1e2waterv}
\end{figure}

The corresponding solution $v$ is shown in 
Fig.~\ref{KBKvgauss1e2waterv}. Note that due to the rescaling of the 
original KBK system (\ref{KBK}), the form of the solitons of 
(\ref{KBKe}) is slightly different. It appears that the oscillations 
traveling to the left in Fig.~\ref{KBKvgauss1e2waterv} are solitons. 

\section{Conclusion}
In this paper we have presented a detailed numerical study of 
solutions to the good  KBK system (\ref{KBK}). The stability of the 
solitons as proven by Angulo \cite{Ang} was illustrated for several 
examples, also for the stationary solution. It was shown that the 
long time behavior of solutions for localized initial data is given 
by solitons plus radiation. In the vicinity of shocks to the 
corresponding dispersionless Saint-Venant system, dispersive shock 
waves were observed. No indication of a blow-up was found even in 
cases where the non-cavitation condition is not satisfied. 

It is an interesting question whether these features can also be 
observed in the 2D variant of the KBK system which is most probably 
not integrable. An important point is whether there are 2D localised 
solitary waves called lumps as in the KP I equation, or whether there 
are no localized 2D structures as for KP II. The 1D solitons 
discussed in the present paper are infinitely extended exact 
solutions to the 2D KBK system, so-called line solitons. There 
stability in the 2D equation has to be studied as well as the 
question whether there can be blow-up in 2D. This will be the subject 
of future work.

 \vspace{0.3cm}

\begin{merci}
We thank Patrik Nabelek for useful discussions on the KBK system.

 The work of the first author  was partially  supported 
by the ANR project ANR-17-EURE-0002 EIPHI. 

Both authors were partially supported by the ANR project ISAAC-ANR-23-CE40-0015-01.
\end{merci}

\bibliographystyle{amsplain}

\end{document}